\documentclass[notitlepage,12pt,a4paper]{article}
\usepackage{booktabs}
\usepackage{subfigure, epsfig, epstopdf}
\usepackage{graphicx,amssymb}
\usepackage{multirow, multicol}
\usepackage{amsmath}
\usepackage{color}
\usepackage{setspace}
\usepackage[left=1in,top=1in,right=1in,bottom=1in, nohead]{geometry}
\usepackage[colorlinks=true,hyperindex,citecolor=blue,linkcolor=blue]{hyperref}
\usepackage{natbib}
\usepackage[utf8]{inputenc}
\usepackage{gensymb}
\usepackage{pdflscape}
\usepackage{courier}
\usepackage{chngpage}
\usepackage{algorithm}
\usepackage{algpseudocode}
\usepackage{longtable}
\usepackage[titletoc]{appendix}
\usepackage{float}
\usepackage{textcomp}
\usepackage{pdflscape}
\usepackage{mathtools}
\usepackage[capitalize]{cleveref}
\usepackage{enumerate}

\allowdisplaybreaks

\newtheorem{proposition}{Proposition}

\newcommand{\qed}{\nobreak \ifvmode \relax \else \ifdim\lastskip<1.5em \hskip-\lastskip \hskip1.5em plus0em minus0.5em \fi \nobreak \vrule height0.75em width0.5em depth0.25em\fi} 

\newcommand{\bzero}{\mathbf 0}

\newcommand{\B}{\mathcal B}

\newcommand{\T}{\mathsf T}

\newcommand{\E}{\mathbb E}

\newcommand{\bmu}{\boldsymbol \mu}

\newcommand{\cov}{\text{Cov}}
\newcommand{\var}{\text{Var}}

\newcommand{\tb}{\text{Beta}}

\newcommand{\X}{\overline{X}}
\newcommand{\Y}{\overline{Y}}
\newcommand{\LX}{\overline{\ln X}}
\newcommand{\XLX}{\overline{X \ln X}}
\newcommand{\XLY}{\overline{X \ln Y}}
\newcommand{\LY}{\overline{\ln Y}}
\newcommand{\YLY}{\overline{Y \ln Y}}

\newcommand{\tgam}{\tilde{\gamma}}

\def\mlx{\overline{\ln X}}
\def\mly{\overline{\ln Y}}

\def\bralpha{\breve{\alpha}}
\def\brbeta{\breve{\beta}}

 \graphicspath{{figure/}} 

\title{Novel Closed-form Point Estimators for the Beta Distribution}

\usepackage{authblk}
\author[1]{Piao CHEN}
\author[2]{Xun XIAO}
\affil[1]{\small{Delft Institute of Applied Mathematics, Delft University of Technology, Netherlands}}
\affil[2]{\small{Department of Mathematics and Statistics, University of Otago, Dunedin, New Zealand}}

\date{}

\begin{document}
\maketitle

\begin{abstract}

In this paper, novel closed-form point estimators of the beta distribution are proposed and investigated. The first estimators are a modified version of Pearson's method of moments. The underlying idea is to involve the sufficient statistics, i.e. log-moments in the moment estimation equations and solve the mixed type of moment equations simultaneously. The second estimators are based on an approximation to Fisher's likelihood principle. The idea is to solve two score equations derived from the log-likelihood function of generalized beta distributions. Both two resulted estimators are in closed-forms, strongly consistent and asymptotically normal. In addition, through extensive simulations, the proposed estimators are shown to perform very close to the ML estimators in both small and large samples, and they significantly outperform the moment estimators. 
\end{abstract}


\section{Introduction}

Statistical point estimation serves as a central role of statistical inference with interval estimation and hypothesisi testing \citep{casella2021statistical}. It dates back to late 19\textsuperscript{th} century when Karl Pearson introduced his method of moment (MM) estimation \citep{pearson1894contributions}. In early 20\textsuperscript{th} century, Ronald A. Fisher suggested the maximum likelihood (ML) estimation  as an alternative to Pearson's MM approach and further established the analytical foundation of MLE and general point estimation in his remarkable paper \citep{fisher1922mathematical}. The pioneer work of these two giants has motivated the development of various estimation approaches in the last century. By and large, the MM estimator is considered to be inferior to the ML estimator since the method of moments are not necessarily based on the sufficient statistics. However, a generalized version of Pearson's idea, i.e. the generalized method of moment (GMM, \citealp{hansen1982large}),  has been advocated by econometricians in light of the restrictive parametric assumptions implied by the likelihood approach. \cite{bera2002mm} made a comprehensive synthesis of different esimation approaches stemming from the original ideas of Pearson and Fisher. 

Although MM estimation and ML estimation yield the same estimators for some commonly-used distributions, e.g., normal distribution, binomial distribution, and Poisson distribution, they usually differ a lot for most distributions. \cite{davidson1974moment} demonstrated that ML estimation is equivalent to a modified version of MM estimation by exploiting the moments of the minimal sufficient statistics for the exponential family. In addition to the intriguing connections between these two classical approach presented in the paper, \cite{davidson1974moment} also mentioned some interesting historical notes on the debate between Pearson and Fisher. 

A major advantage of MM estimation is that it tends to be less computatonally involved than ML estimation. For example, obtaining ML estimates for gamma distribution and beta distribution relies heavily on iterative numerical algorithms to maximize the likelihood function, while MM estimates enjoys some simple closed forms. Due to the rapid development of high-speed digital techniques in the last century, the drawback of ML estimation in computational efficiency has been greatly relieved in most scenarios. However, a massive amout of data comes into the scope in a more and more complicated fashion. The numerical algorithms for finding ML estimates can become hindrance in certain cases. \citet{huang2019distributed} studied a one-step estimator as a surrogate to the ML estimator which needs multiple step iterations under the framework of distribute computing. \cite{xiao2021computing} showed that maximizing the log-likehiood in change point analysis can be rather time-consuming. \cite{xiong2021faster} argued that using closed-form estimator can accelerate the inference of gamma mixture model.

Therefore, some recent papers have proposed and investigated new closed-form estimators for gamma distribution, beta distribution and other distributions. \cite{ye2017closed} firstly reported new closed-from estimators for gamma distribution and showed that the performance of new estimators are very close to ML estimators. Many follow-up works have been published since then. Similar new closed-form estimators have also been found for matrix-variate Gamma distribution \citep{alfelt2020closed}, weighted Lindley distribution \citep{kim2021new}, and Nakagami distribution \citep{zhao2021closed}. 
 \cite{tamae2020score} demonstrated that Ye-Chen estimator for gamma distribution can be derived by following the idea of approximating the special functions in score equations via mixed moments.
 In light of the close connection between beta distribution and gamma distribution, they further applied their idea to beta distribution to obtain some interesting closed-form estimators called score-adjusted estimators and score-adjusted moment estimators. 
\cite{papadatos2022point} constructed the unbiased version of Ye-Chen estimators for gamma distribution and then derived the closed-form estimators for beta distribution using Stein's identity and $U$-statistics. Surprisingly, the estimators presented by \cite{tamae2020score} agree with ones presented by \cite{papadatos2022point} for beta distribution. 

Although both \cite{tamae2020score} and \cite{papadatos2022point} found the new closed-form estimators for beta distribution, their technical details are considerably different to the original approach adopted by \cite{ye2017closed}. In this paper, we will restrict our interest to deriving novel closed-form estimators for beta distribution by tracing the closed-form estimators for the gamma distribution proposed in \cite{ye2017closed}. 

Mathematically, the beta distribution is indexed by two parameters $\alpha$ and $\beta$, and the probability density function (PDF) of $\tb(\alpha,\beta)$ is given by 
\begin{equation}
\label{eq:betapdf}
f(x;\alpha,\beta) = \frac{x^{\alpha-1}(1-x)^{\beta-1}}{\B(\alpha,\beta)}, 0<x<1, 
\end{equation}
where $\alpha>0$, $\beta>0$, $\B(\alpha,\beta)$ is the beta function. As seen from the PDF, the beta distribution has the property known as the mirror-image symmetry: if $X\sim \tb(\alpha,\beta)$, then $Y \equiv 1-X \sim \tb(\beta,\alpha)$. 

In the literature, estimating the beta distribution is often based on two classical methods, i.e., Fisher's ML estimation and Pearson's MM estimation. Let $X_1,\dots,X_n$ be $n$ i.i.d. copies of $X$, where $X \sim \tb(\alpha,\beta)$. In addition, let $Y_i = 1-X_i, i=1,\dots,n$. The ML estimators of $\alpha$ and $\beta$ are given by the solution of the following system of score equations
\begin{equation}
\label{eq:mle}
\begin{aligned}
n[\overline{\ln X}+\psi(\alpha+\beta)-\psi(\alpha) ] &= 0 \\ 
n[\overline{\ln Y}+\psi(\alpha+\beta)-\psi(\beta) ] &=0,
\end{aligned}
\end{equation}
where $\overline{\cdot}$ denotes the sample mean of the corresponding random sample and $\psi (z)=\frac{d}{dz} \ln \Gamma(z) $ is the digamma function. As seen, the ML estimators of the beta parameters do not have closed-forms, and they have to be obtained by iterative algorithms. As noted by many authors \citep[e.g.,][]{lau1991effective}, starting values are crucial for the efficient convergence of most iterative algorithms, which is generally not an easy task. Moreover, even with good starting values, most numerical algorithms fail to converge when $\exp(\overline{\ln X}) + \exp(\overline{\ln Y}) > 0.95$ \citep[see][]{cordeiro1997bias}, which significantly impedes the use of the ML estimation. Lastly, the number of iterations required in obtaining the ML estimators can be quite large in some cases, which somehow hinders the application of the beta distribution when fast computing is required; see \citet{huang2019distributed} for an example on estimating the beta distribution in a large distributed system.

On the other hand, the MM estimation makes use of the following two moment equations
\begin{equation}
\label{eq:momeq}
\E[X] = \frac{\alpha}{\alpha+\beta} \quad \text{and} \quad \E[X^2] = \frac{\alpha(\alpha+1)}{(\alpha+\beta)(\alpha+\beta+1)}.
\end{equation}
By using the sample moments $\overline{ X}$ and $\overline{X^2}$ in place of $\E[X]$ and $\E(X^2)$, the moment estimators are then obtained as
\begin{equation}
\label{eq:mom}
\hat \alpha_{\text{mom}} = \frac{\overline{X}(\overline{X} - \overline{ X^2})}{\overline{ X^2}-\overline{ X}^2} \quad \text{and} \quad
\hat \beta_{\text{mom}}  = \frac{\overline{Y}(\overline{Y} - \overline{ Y^2})}{\overline{ Y^2}-\overline{ Y}^2}.
\end{equation}
As one may expect, the moment estimators are not efficient under either small or large samples.

In view of the deficiencies of the classical methods, this study aims to propose more efficient versions of closed-form estimators for the beta distribution. Particularly, \cite{ye2017closed} derived their closed-form estimators for the gamma distribution based on two out of three score equations of generalized gamma distribution. It is noticed that Ye-Chen estimators can be interpreted from the perspective of either Pearson or Fisher with certain modifications. 
So, we firstly offer a Pearson type interpretation of Ye-Chen estimators as mixed moment estimators involving the first order moment equation and an equation related to the covariance between the moment and the log-moment of the gamma distribution. By following this interpretation, we replace the second moment equation in \eqref{eq:momeq} by an equation related to the covariance between the moments and the log-moments of the beta distribution. Solving these two mixed moment equations leads to the Pearson type of closed-form estimators for beta distribution, i.e. the mixed moment estimators. 

Interestingly, these mixed moment estimators coincide with the score-adjusted moment estimators reported in \cite{tamae2020score} and \cite{papadatos2022point}. However, both mixed moment estimators and score-adjusted moment estimators construct the closed-form estimators in a somewhat ad hoc manner. The related details will be further examined in further discussions. With the aid of the mirror-image symmetry property of the beta distribution, a more elegant type of closed-form estimators for the beta distribution is derived from the score equations of generalized beta distribution just like the work in \cite{ye2017closed}. Moreover, it is shown that the second type of closed-form estimators are a refined version of the ill-posed score-adjusted estimators in \cite{tamae2020score}. Therefore, the second type of closed-form estimators are essentially refined score-adjusted estimators. In addition, we show that both two types of proposed estimators are well-defined, strongly consistent, and asymptotically normal distributed. Numerical results suggest that the asymptotic variances of both two types of new estimators are very close to the Cram\'{e}r-Rao lower bound. Monte Carlo simulation studies are also conducted to assess the performance of the proposed closed-form estimators under both small and large samples.

The rest of the paper is organized as follows.
\cref{sec:idea} first briefly reviews the ideas of \cite{ye2017closed} and \cite{tamae2020score} for the purposes of motivation and comparison. It further presents an alternative derivation of the closed-form SAM estimators for beta distributions with their asymptotic properties. 
\cref{sec:newest} proposes the new closed-form estimators for beta distributions and investigates their properties. 
\cref{sec:numstudy} conducts numerical studies to assess the performance of proposed estimators. At last, \cref{sec:conclusion} concludes the paper.

\section{Moment-type estimators for beta distribution}
\label{sec:idea}

\cite{ye2017closed} derived closed-form estimators for the gamma distribution with shape parameter $k>0$ and scale parameter $\theta>0$, denoted by $\text{Gamma}(k,\theta)$, by utilising two out of three score equations of generalized gamma distribution. \cite{tamae2020score} presented a score-adjusted interpretation of the closed-form gamma estimators in \cite{ye2017closed} and generalized their score-adjusted idea to the beta distribution. Their core idea is to replace digamma functions, i.e., the main hindrance in solving the score equations \eqref{eq:mle}, by the sample log-moments of the beta distribution. However, they found that the corresponding score-adjusted estimators for the beta distribution show poor performance when $\alpha\approx\beta$ and become even non-identifiable if $\alpha=\beta$. Therefore, \cite{tamae2020score} constructed the closed-form score-adjusted moment (SAM) estimators by combining their score-adjusted estimators with the first order moment condition in \cref{eq:momeq}. 
The SAM estimators are also derived in \cite{papadatos2022point} by using the well-known Stein-type identity. 

Given a random sample from $\text{Gamma}(k,\theta)$, \cite{xiao2021computing} noticed that the Ye-Chen estimators can be written as 
\begin{align}
\hat{k} &= \overline{X}/\hat{\theta}, \label{eq:gammamoment}\\ 
\hat{\theta} &= \overline{X\ln X}-\overline{X}\cdot\overline{\ln X}. \label{eq:gammalogmoment}
\end{align}
Here, \cref{eq:gammamoment} is simply the first order moment condition. \cref{eq:gammalogmoment} is a biased estimator of $\text{Cov}(X,\ln X)=\theta$ which can be regarded as a mixed moment condition based on the first order moment and the first order log moment. 

Interestingly, similar moment conditions can also be found in the beta distribution. Let $X\sim \tb(\alpha,\beta)$. Recall that $Y \equiv 1-X \sim \tb(\beta,\alpha)$ by the mirror-image symmetry of the beta distribution.  Let $X_1,\dots,X_n$ be $n$ i.i.d. copies of $X$ and $Y_i = 1-X_i, i=1,\dots,n$. Our proposed estimators are motivated by the following fact
\begin{align*}
\text{Cov}(X, \ln X)  &=\mathbb{E}[X \ln X]-\mathbb{E} [X] \mathbb{E} [\ln X] = \frac{\beta}{(\alpha+\beta)^{2}},  \\ 
\text{Cov}(Y, \ln Y)  &=\mathbb{E}[Y \ln Y]-\mathbb{E} [Y] \mathbb{E} [\ln Y] = \frac{\alpha}{(\alpha+\beta)^{2}}.
\end{align*}
By adding up the above two equations, we observe that 
\begin{equation}
\label{eq:logmoment}
\mathbb{E}[X \ln X]-\mathbb{E} [X] \mathbb{E} [\ln X]+\mathbb{E}[Y \ln Y]-\mathbb{E} [Y] \mathbb{E} [\ln Y] = \frac{1}{\alpha+\beta} .
\end{equation}
Using \eqref{eq:logmoment} in place of the second moment equation in \eqref{eq:momeq}, the new estimation equations become
\begin{equation}
\label{eq:esteq}
      \left\{
                \begin{array}{ll}
                  \E[X]={\alpha}/{(\alpha+\beta)}\\
\mathbb{E}[X \ln X]-\mathbb{E} [X] \mathbb{E} [\ln X]+\mathbb{E}[Y \ln Y]-\mathbb{E} [Y] \mathbb{E} [\ln Y] = {1}/{(\alpha+\beta)}
                \end{array}
              \right.
\end{equation}
Replacing the expectations by the corresponding sample means, closed-form estimators for beta distribution are given by
\begin{align} 
\tilde{\alpha} &=\frac{\overline{X}}{\overline{X \ln X}-\overline{X} \cdot \overline{\ln X}+\overline{Y \ln Y}-\overline{Y} \cdot \overline{\ln Y}} , \label{eq:estalpha}\\ 
\tilde{\beta} &=\frac{\overline{Y}}{\overline{X \ln X}-\overline{X} \cdot \overline{\ln X}+\overline{Y \ln Y}-\overline{Y} \cdot \overline{\ln Y}} . \label{eq:estbeta}
\end{align}

As seen, the new estimators involve the sample means of $X$, $\ln X$, $\ln Y$, $X \ln X$ and $Y \ln Y$, and hence they can be treated as the mixed type of moment estimators. These estimators for the beta distribution agree with the closed-form SAM estimators presented in \cite{tamae2020score} and \cite{papadatos2022point}. \cite{tamae2020score} offered an intuitive explanation of the potential nice performance of these closed-form estimators. Particularly, the estimating equation in \eqref{eq:logmoment} is the difference of two score-adjusted equations studied in \cite{tamae2020score}. For gamma distribution, the score-adjusted equations can be regarded as a stochastic approximation to the original score equations and the approximation error can be very small \citep{xiao2021computing}. The first order moment condition is then added to obtain a solvable system of equations. We shall notice that, though the sample mean $\bar{X}$ is not an asymptotic efficient estimator of the population mean $\E[X]$ for the beta distribution, it is expected to be very robust as an unbiased estimator especially for small sample sizes. In other words, the first order moment condition in \eqref{eq:esteq} may borrow robustness to the SAM estimators under small sample sizes while the second estimating equation based on log-moments in \eqref{eq:esteq} constrains the closed-form estimators in a neighbourhood of the ML estimators. 

In addition to the desirable closed-forms, these estimators have the following standard properties summarised in \cref{Proposition1}. The technical proofs are given in \cref{proof:thm1}. Particularly,  the covariance matrix of these estimators is also given in \cite{tamae2020score} but it is derived in a different approach in this paper.

\begin{proposition} \label{Proposition1}
\hfill
\begin{enumerate}[i.]
\item The proposed estimators $\tilde \alpha$ and $\tilde \beta$ in \eqref{eq:estalpha} and \eqref{eq:estbeta} are well-defined and positive if there exist at least two distinct observations in a random sample with sample size $n\ge 2$ from the unit interval $(0,1)$.
\item The proposed estimators $\tilde \alpha$ and $\tilde \beta$ in \eqref{eq:estalpha} and \eqref{eq:estbeta} are strongly consistent.
\item The proposed estimators $\tilde \alpha$ and $\tilde \beta$ in \eqref{eq:estalpha} and \eqref{eq:estbeta} are asymptotically normal distributed as $n \to \infty$. In specific, 
\begin{equation*}
\sqrt{n} \left[\binom{\tilde \alpha}{\tilde \beta}-\binom{\alpha}{\beta}\right] \to^{d} \mathsf{MVN}(\boldsymbol{0},\Sigma_1),
\end{equation*}
where the covariance matrix $\Sigma_1$ is 
\footnotesize
\begin{equation*}
\frac{1}{\alpha+\beta+1}\left[\begin{array}{ll}{\alpha^{3} \beta[\psi_1(\alpha)+ \psi_1(\beta)]+\alpha^{2}(\alpha+\beta+1)-\alpha \beta}, & 
{(\beta-1) \alpha^{2}+\alpha^{2} \beta^{2}[\psi_1(\alpha)+ \psi_1(\beta)]+(\alpha-1)\beta^{2}} \\ 
{(\beta-1) \alpha^{2}+\alpha^{2} \beta^{2}[\psi_1(\alpha)+ \psi_1(\beta)]+(\alpha-1)\beta^{2}},  &
 {\alpha \beta^{3}[\psi_1(\alpha)+ \psi_1(\beta)]+\beta^{2}(\alpha+\beta+1)-\alpha \beta}\end{array}\right]
\end{equation*}
\normalsize
and $\psi_1 (z)=\frac{d^2}{dz^2} \ln \Gamma(z) $ is the trigamma function.
\end{enumerate}
\end{proposition}




\section{Refined score-adjusted estimators}
\label{sec:newest}

Many researchers \citep{tamae2020score,xiao2021computing,papadatos2022point} have provided interesting and meaningful interpretations on Ye-Chen estimator and further attempted to generalize their ideas to the case of the beta distribution. Eventually, these efforts converge to the SAM estimators proposed by \cite{tamae2020score} as discussed in \cref{sec:idea}. However, the original idea of \cite{ye2017closed}, i.e., the intriguing connections between the gamma distribution and the generalized gamma distribution, has never been visited. 

Just like the gamma distribution is a special case of the generalized gamma distribution, the beta distribution can be regarded as a special case of the generalized beta distribution with the PDF
\[
f_{\text{gb}}(x;\alpha,\beta,r)=\frac{rx^{r\alpha-1}(1-x^r)^{\beta-1}}{\B(\alpha,\beta)}, x\in(0,1).
\]
by setting $r=1$. In other words, $X\sim\text{Beta}(\alpha,\beta)$, $X^r\sim\text{GB}(\alpha,\beta,r)$. 

Borrowing the idea of \cite{ye2017closed}, we derive the score equations of generalized beta distribution and then set $r=1$ as follows:
\begin{align}
n[\mlx+\psi(\alpha+\beta)-\psi(\alpha)]&=0, \label{eq:sb1}\\
n[\mly+\psi(\alpha+\beta)-\psi(\beta)]&=0, \label{eq:sb2}\\
n[1+\alpha\mlx-(\beta-1)\overline{X\ln X/(1-X)}]&=0. \label{eq:sb3}
\end{align}

Unfortunately, the idea of \cite{ye2017closed} does not work directly since the first two equations \cref{eq:sb1,eq:sb3} are just the beta score equations in \cref{eq:mle} which involve the digamma function. The last equation \cref{eq:sb3} is nice to handle but two out of three equations fail to yield any closed-form estimator. One may add the first order moment condition $\bar{X}=\alpha/(\alpha+\beta)$ and then solve two equations together. But the resulted estimator looks highly unreliable in the absence of the other sufficient statistics $\overline{\ln Y}$. In light of integrating $\overline{\ln Y}$ into the estimation procedure, this flaw can be remedied by recalling the  the property of mirror-image symmetry. Making the transformation $Y=1-X$, we can get another independent equation from the score equations of generalized beta distribution for $Y^r\sim\text{GB}(\beta,\alpha,r)$ as follows
\begin{equation}
n[1+\beta\mly-(\alpha-1)\overline{Y\ln Y/(1-Y)}]=0. \label{eq:sb4}
\end{equation}

Solving \cref{eq:sb3,eq:sb4} yields new closed-form estimators for beta distribution as:
\begin{align}
\bralpha&=\frac{(1+M_X)\LY+(1+M_Y)M_X}{M_XM_Y-\LX\cdot\LY}, \label{eq:estalpha2}\\
\brbeta&=\frac{(1+M_Y)\LX+(1+M_X)M_Y}{M_XM_Y-\LX\cdot\LY}. \label{eq:estbeta2}
\end{align}
where $M_X=\overline{X\ln X/(1-X)}$ and $M_Y=\overline{Y\ln Y/(1-Y)}$.

Similarly to the SAM estimators in \cref{eq:estalpha,eq:estbeta}, the new estimators involve the sample means of $\ln X$, $\ln Y$, $X \ln X/(1-X)$ and $Y \ln Y/(1-Y)$, and hence they can also be treated as the mixed type of moment estimators. In addition to the desirable closed-forms, the new estimators have the following properties summarised in \cref{Proposition2}. Rather than the standard delta method used in proving, a different approach based on generalized method of moment is used to find the covariance matrix $\Sigma_2$. The technical proofs are presented in \cref{proof:thm2}. We must admit that the form of the covariance matrix $\Sigma_2$ is a little bit too cumbersome. However, further simplifying the whole expression is highly nontrivial with digamma and trigamma functions involved.

\begin{proposition} \label{Proposition2}
\hfill
\begin{enumerate}[(i)]

\item The proposed estimators $\bralpha$ and $\brbeta$ in \eqref{eq:estalpha2} and \eqref{eq:estbeta2} are well-defined and positive if there exist at least two distinct observations in a random sample with sample size $n\ge 2$ from the unit interval $(0,1)$.

\item The proposed estimators $\bralpha$ and $\brbeta$ in \eqref{eq:estalpha2} and \eqref{eq:estbeta2}  are strongly consistent.

\item The proposed estimators $\bralpha$ and $\brbeta$ in \eqref{eq:estalpha2} and \eqref{eq:estbeta2}  are asymptotically normal distributed as $n \to \infty$. In specific, 
\begin{equation*}
\sqrt{n} \left[\binom{\bralpha}{\brbeta}-\binom{\alpha}{\beta}\right] \to^{d} \mathsf{MVN}(\boldsymbol{0},\boldsymbol \Sigma_2),
\end{equation*}
where the covariance matrix is given by 
\begin{align*}
\boldsymbol \Sigma_2&=\frac{1}{\left(\kappa_{\alpha\beta}\kappa_{\beta\alpha}-\tau_{\alpha\beta}\tau_{\beta\alpha}\right)^2}
\begin{bmatrix} 
\kappa_{\beta\alpha} & \tau_{\alpha\beta}\\ 
\tau_{\beta\alpha} & \kappa_{\alpha\beta}
\end{bmatrix}
\begin{bmatrix} 
\omega_{\alpha\beta} & \rho \\ 
\rho & \omega_{\beta\alpha}
\end{bmatrix}
\begin{bmatrix} 
\kappa_{\beta\alpha} & \tau_{\beta\alpha}\\ 
\tau_{\alpha\beta} & \kappa_{\alpha\beta}
\end{bmatrix}
\end{align*}
with
\begin{align*}
\kappa_{\alpha\beta}&=\psi(\alpha)-\psi(\alpha+\beta),\\ 
\tau_{\alpha\beta}&=\frac{\alpha}{\beta-1}[\psi(\alpha+1)-\psi(\alpha+\beta)],\\
\omega_{\alpha\beta}&=\frac{\beta}{\beta-2}+\frac{\alpha(\alpha+\beta-1)}{\beta-2}[\psi_1(\alpha)-\psi_1(\alpha+\beta)+\kappa_{\alpha\beta}^2]+\frac{2(2\alpha+\beta-1)}{\beta-2}\kappa_{\alpha\beta},\\
\rho&=1+\tau_{\alpha\beta}+\tau_{\beta\alpha}+(\alpha+\beta+1)[\psi_1(\alpha+\beta)-\kappa_{\alpha\beta}\kappa_{\beta\alpha}].
\end{align*}

\end{enumerate}
\end{proposition}

Starting from two score equations of the beta distribution in \cref{eq:mle}, \cite{tamae2020score} replaced two differences of digamma functions by two mixed moments as
\begin{align}
\E[\ln X]=\psi(\alpha)-\psi(\alpha+\beta)&=\frac{(\alpha+\beta)\E[X\ln X]}{\alpha}+\frac{1}{\alpha+\beta}-\frac{1}{\alpha}, \label{eq:sa1}\\
\E[\ln Y]=\psi(\beta)-\psi(\alpha+\beta)&=\frac{(\alpha+\beta)\E[X\ln Y]}{\alpha}+\frac{1}{\alpha+\beta}. \label{eq:sa2}
\end{align}
Without involving the first order moment condition, they obtained the so-called score-adjusted estimators for the beta distribution by plugging in the sample counterparts of the moments and solving \cref{eq:sa1,eq:sa2}. However, they further noticed that the score-adjusted estimators perform poorly when $\alpha$ is close to $\beta$ since \cref{eq:sa1,eq:sa2} are identical if $\alpha=\beta$. To remedy the identifiability of score-adjusted estimators, they combined \cref{eq:sa1,eq:sa2} with the first order moment condition and obtained the final SAM estimators in a somewhat ad hoc way. 

Interestingly, although the score-adjusted method fails to deliver meaningful closed-form estimators for the beta distribution, we notice that our new closed estimators motivated by Ye-Chen estimators can be regarded as a refined version of score-adjusted estimators. Particularly, \cref{eq:sb3,eq:sb4} can be rewritten as 
\begin{align}
\mlx &= \frac{(\beta-1)\overline{X\ln X/(1-X)}}{\alpha}-1, \label{eq:sars1}\\
\mly &= \frac{(\alpha-1)\overline{Y\ln Y/(1-Y)}}{\beta}-1. \label{eq:sars2}
\end{align}
It is not hard to verify that 
\begin{align}
\E[\ln X] &=\psi(\alpha)-\psi(\alpha+\beta) = \frac{(\beta-1)\E[X\ln X/(1-X)]}{\alpha}-\frac{1}{\alpha}, \label{eq:sar1}\\
\E[\ln Y] &=\psi(\beta)-\psi(\alpha+\beta) = \frac{(\alpha-1)\E[Y\ln Y/(1-Y)]}{\beta}-\frac{1}{\beta}, \label{eq:sar2}
\end{align}
by noticing that 
\begin{align}
\E[X\ln X/(1-X)]&=\frac{\alpha}{\beta-1}[\psi(\alpha+1)-\psi(\alpha+\beta)], \label{eq:betamoment1}\\
\E[Y\ln Y/(1-Y)]&=\frac{\beta}{\alpha-1}[\psi(\beta+1)-\psi(\alpha+\beta)]. \label{eq:betamoment2}
\end{align}

\cref{eq:sar1,eq:sar2} indicate that our estimators can be regarded as another version of score-adjusted estimators. The score-adjusted equation is a stochastic approximation to the original score equations driven by Fisher's likelihood principle. Hence, the resulted estimators are naturally an approximation to the ML estimator.

Since \cref{eq:sb3,eq:sb4} are both linear equations of $\alpha$ and $\beta$, it is not hard to verify that the system of linear equations is determined. A proof has been provided in \cref{Proposition2} (i).  One shall also notice that, by L'H\^{o}pital's rule, two expectations in \cref{eq:betamoment1,eq:betamoment2} are still well-defined if $\alpha$ or $\beta$ approaches one as 
\begin{align}
\E[X\ln X/(1-X)]&=-\alpha\psi_1(\alpha+1),~~~\beta=1,\\
\E[Y\ln Y/(1-Y)]&=-\beta\psi_1(\beta+1),~~~\alpha=1.
\end{align} Therefore, this version of score-adjusted estimators solves the problem of identifiability when $\alpha=\beta$ in the original version of score-adjusted estimators proposed by \cite{tamae2020score}.

\section{Numerical study}
\label{sec:numstudy}

\subsection{Numerical evaluation}

We compare the asymptotic variance of the two proposed estimators, the ML estimators and the moment estimators numerically. The asymptotic covariance of the ML estimators in \eqref{eq:mle}, also known as the Cram\'{e}r-Rao  bound, is given by 
\begin{equation*}
\frac{1}{\psi_1(\alpha)\psi_1(\beta)-[\psi_1(\alpha)+\psi_1(\beta)]\psi_1(\alpha+\beta)}
\left[\begin{array}{cc}
\psi_1(\beta)-\psi_1(\alpha+\beta) & 
\psi_1(\alpha+\beta) \\ 
\psi_1(\alpha+\beta) &
\psi_1(\alpha)-\psi_1(\alpha+\beta)
\end{array}\right],
\end{equation*}
where $\psi (z)=\frac{d}{dz} \ln \Gamma(z) $ is the digamma function and $\psi_1 (z)=\frac{d^2}{dz^2} \ln \Gamma(z) $ is the trigamma function. On the other hand,  the asymptotic covariance  of the moment estimators in \eqref{eq:mom} is given by
\small
\begin{equation*}
\frac{(\alpha+1) (\beta+1) (\alpha+\beta)}{(\alpha+\beta+1)(\alpha+\beta+2)(\alpha+\beta+3)}
\left[\begin{array}{cc}
M(\alpha,\beta) & 
1+\alpha+\beta+2(\alpha+\beta)^2 \\ 
1+\alpha+\beta+2(\alpha+\beta)^2 &
M(\beta,\alpha)
\end{array}\right],
\end{equation*}
\normalsize
where $M(\alpha,\beta) ={\alpha}[1+(2\beta+3)\alpha^2+(2\beta^2+4\beta+5)\beta+(4\beta^2+7\beta+4)\alpha]/[{\beta(\beta+1)}]$. 

\begin{figure}
\centering
\subfigure{\includegraphics[width=5.2cm]{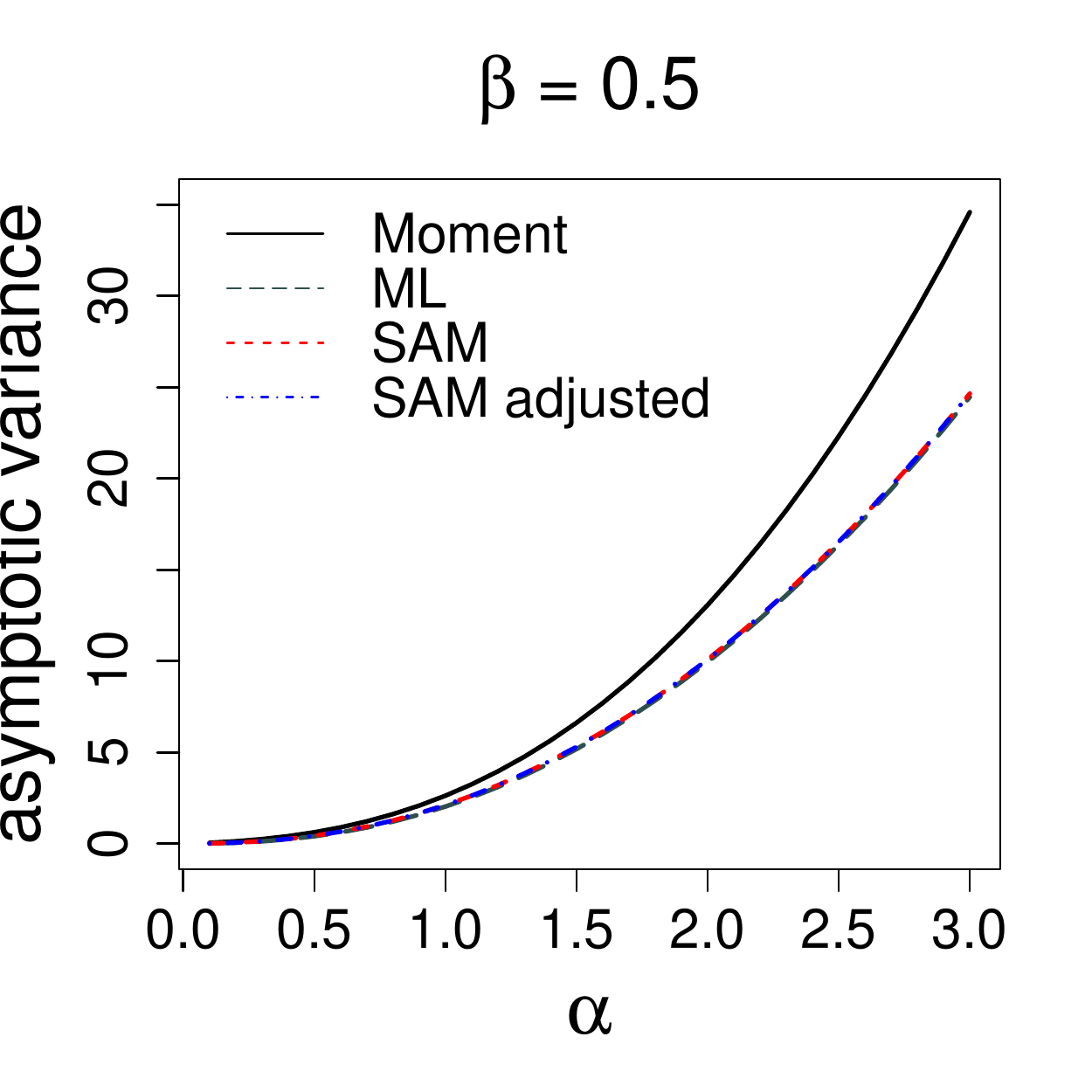}}
\subfigure{\includegraphics[width=5.2cm]{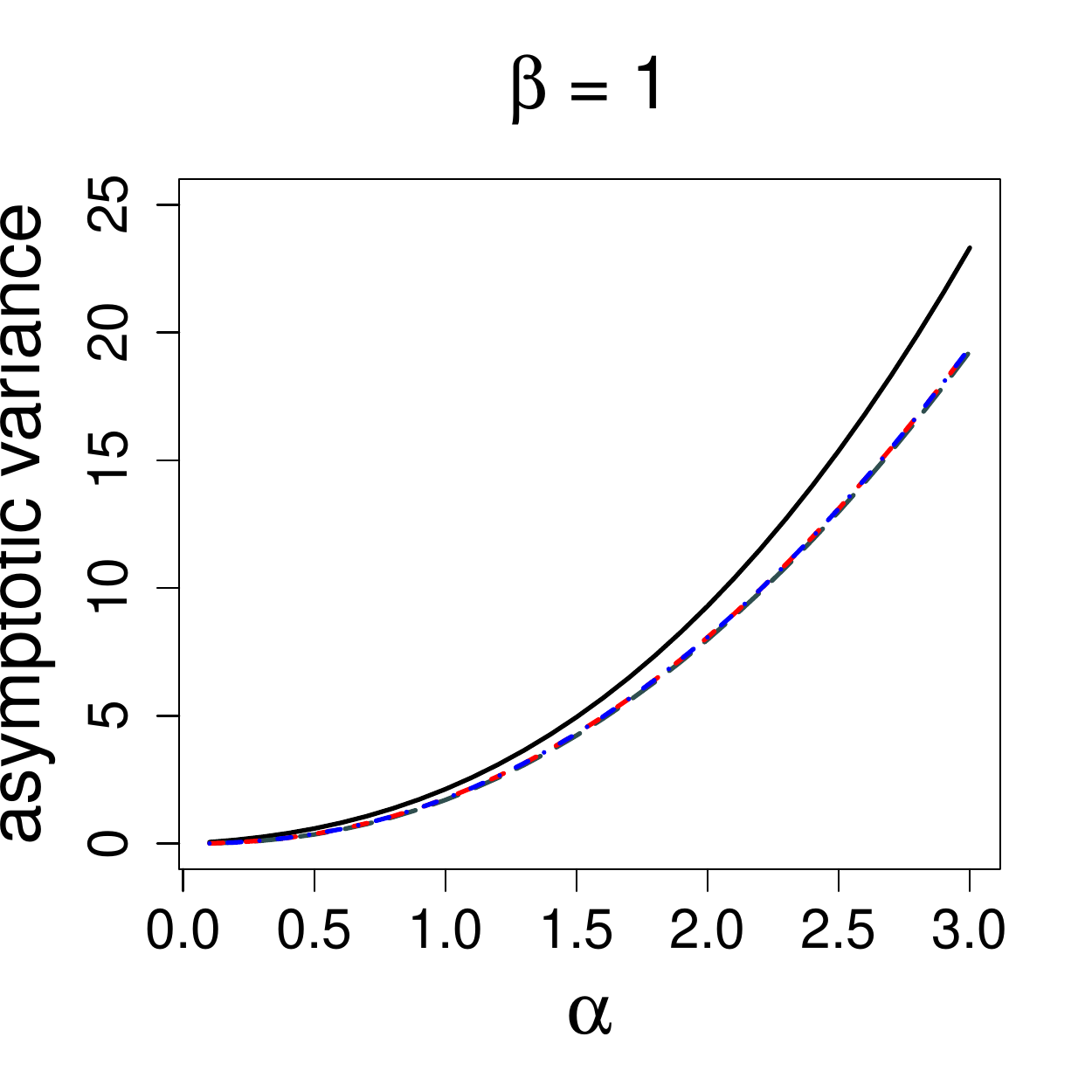}}
\subfigure{\includegraphics[width=5.2cm]{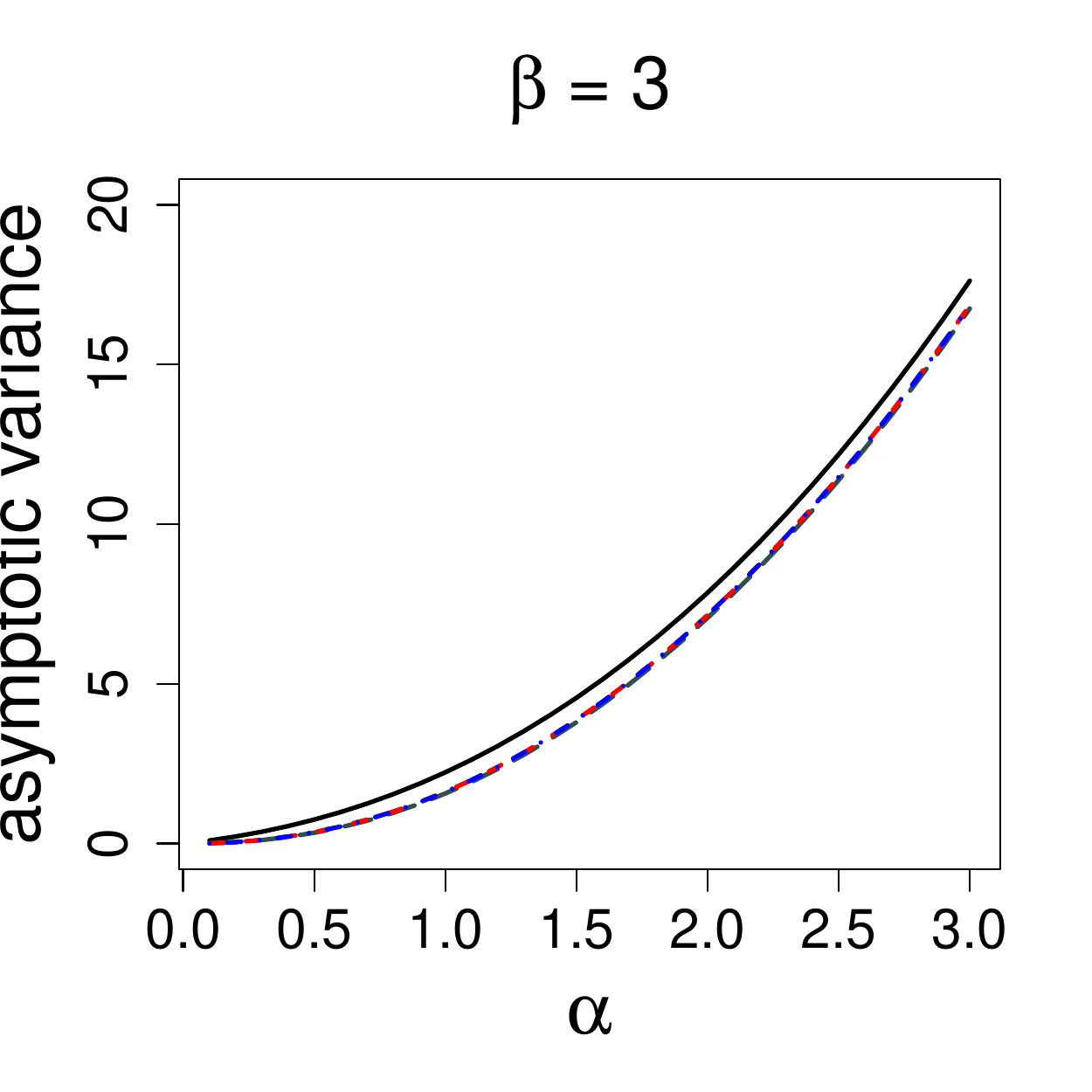}}
\caption{Asymptotic variance of the moment estimator, the ML estimator and the two proposed estimators in estimating $\alpha$ when $\beta=0.5,1,3$.}
\label{fig:avar}
\end{figure}

Due to the mirror-image symmetry, we only consider the asymptotic variance of estimators of $\alpha$ under $\beta=0.5,1,3$. The graphical results are shown in Figure~\ref{fig:avar} with the true $\alpha$ varying from 0.1 to 3.  As seen, although the proposed estimators have slightly larger asymptotic variance than the ML estimator,  they are almost indistinguishable in the figures. In addition, 
both the proposed estimators and the ML estimator  outperform the moment estimator by a large margin. 

The numerical results indicate that the proposed estimators are nearly as efficient as the ML estimators asymptotically. 
Nevertheless, the closed-form expressions make the proposed estimators preferable in practice as iteration algorithm is required to obtain the ML estimators, which can be computationally expensive in some applications. For example, \cite{xiao2021computing} demonstrated that the computation efficiency of detecting the multiple change points within successive observations can be greatly improved by plugging a closed-form estimator into the log-likelihood function for the gamma distribution. 
If multiple change points detection of the beta distribution is considered, the proposed closed-form estimators should also be very useful in improving the computation efficiency.

\subsection{Simulation Study}

\begin{figure}
\centering
\subfigure{\includegraphics[width=7cm]{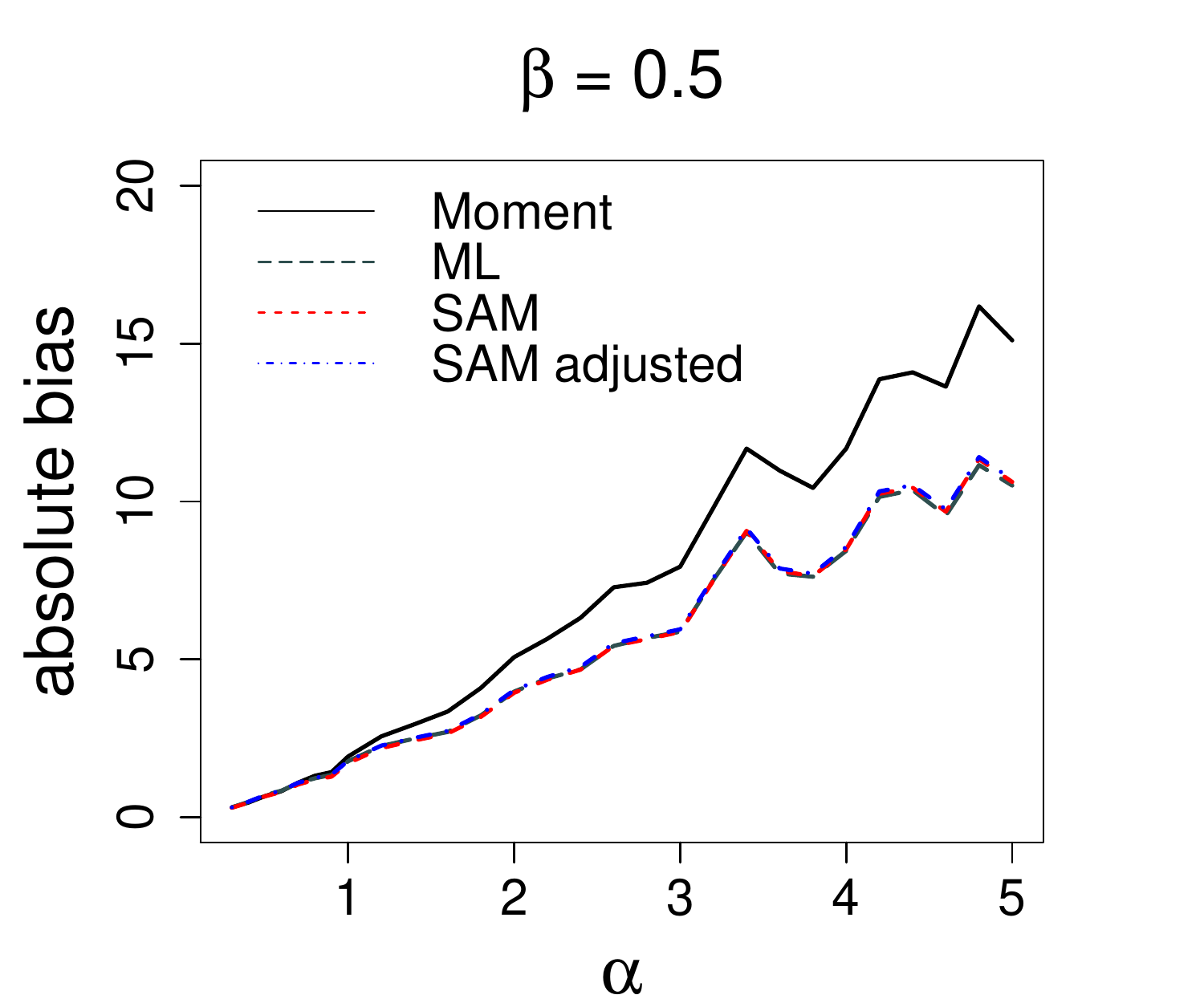}}
\subfigure{\includegraphics[width=7cm]{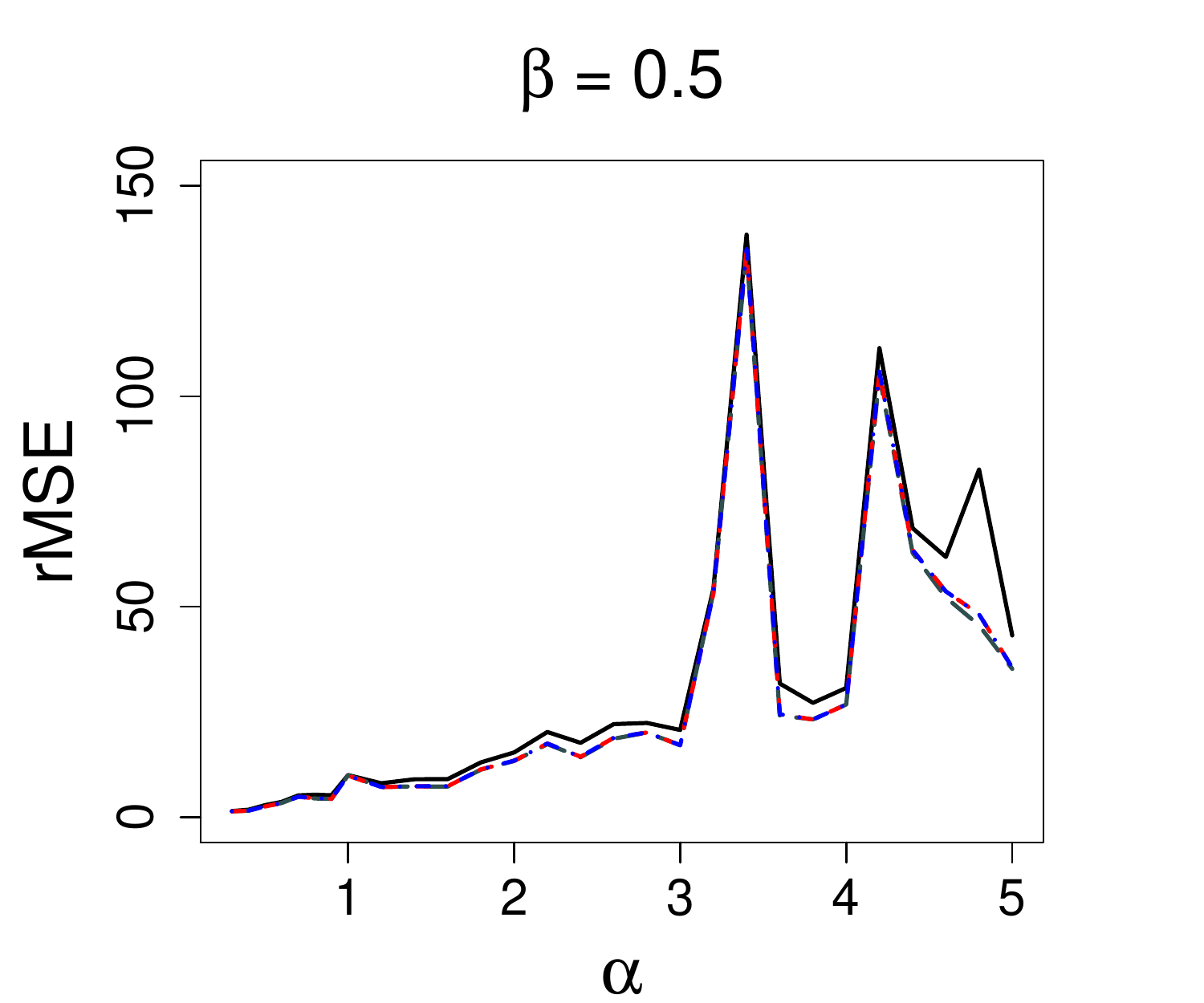}}
\hfill
\\
\subfigure{\includegraphics[width=7cm]{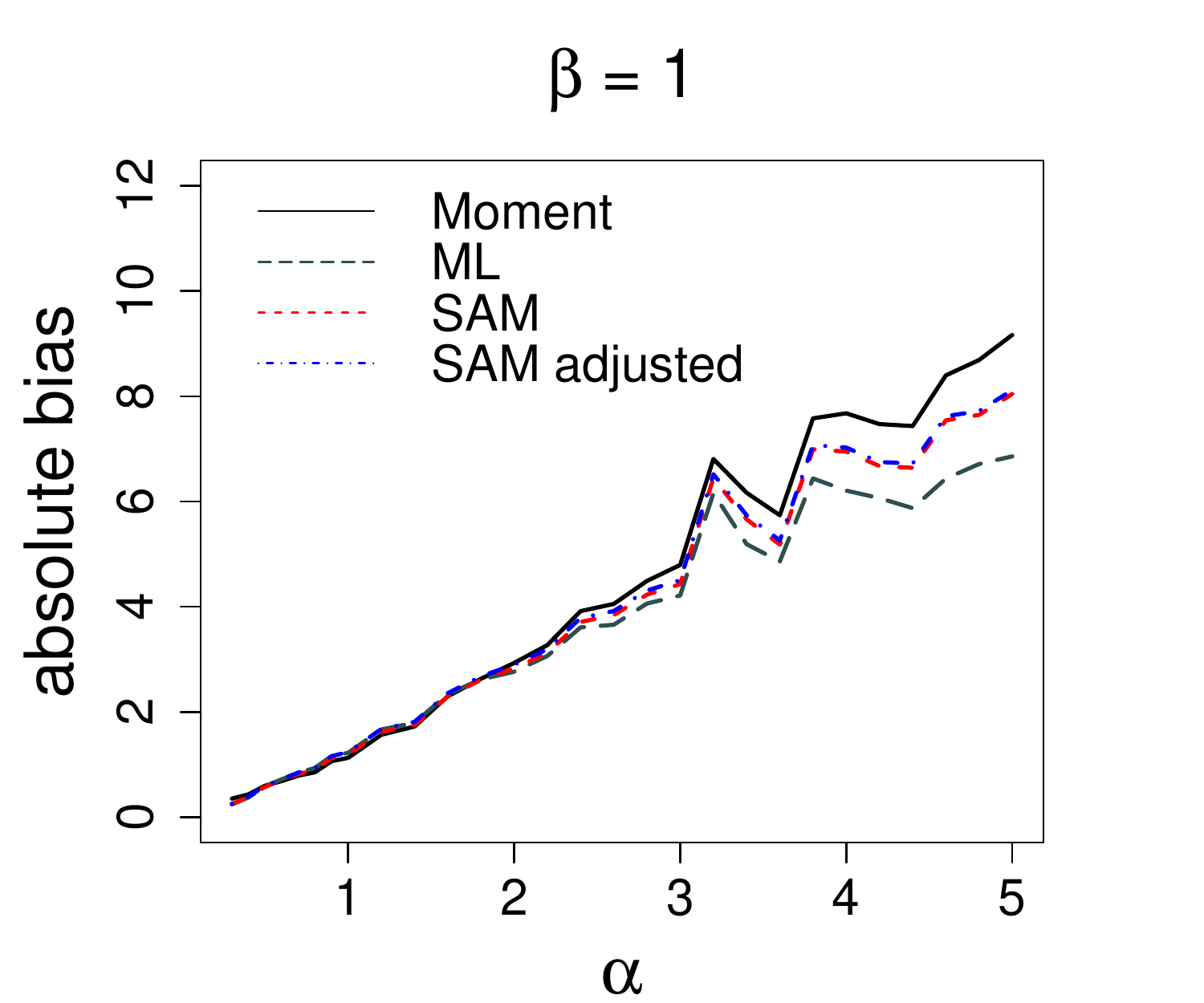}}
\subfigure{\includegraphics[width=7cm]{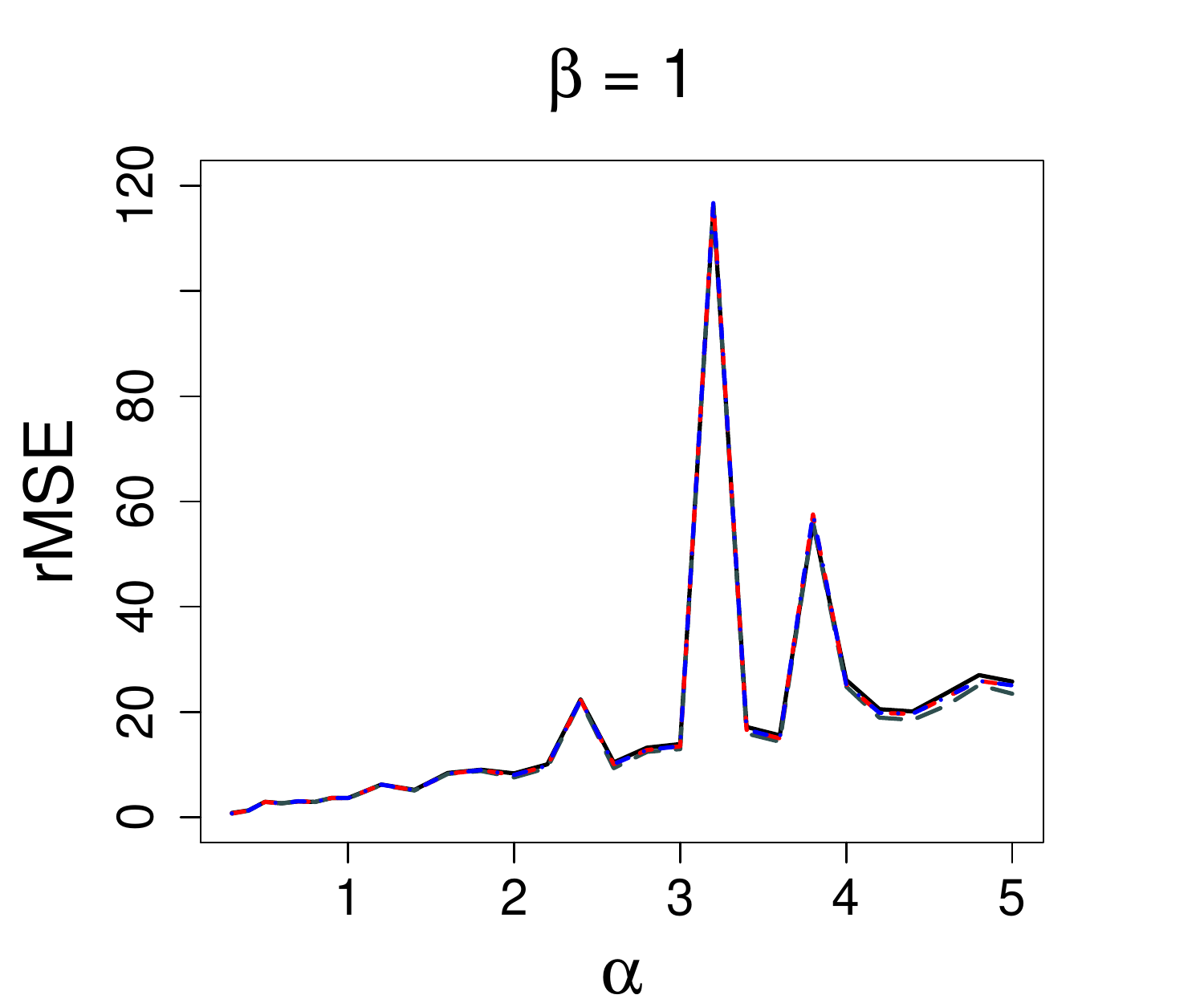}}
\hfill
\\
\subfigure{\includegraphics[width=7cm]{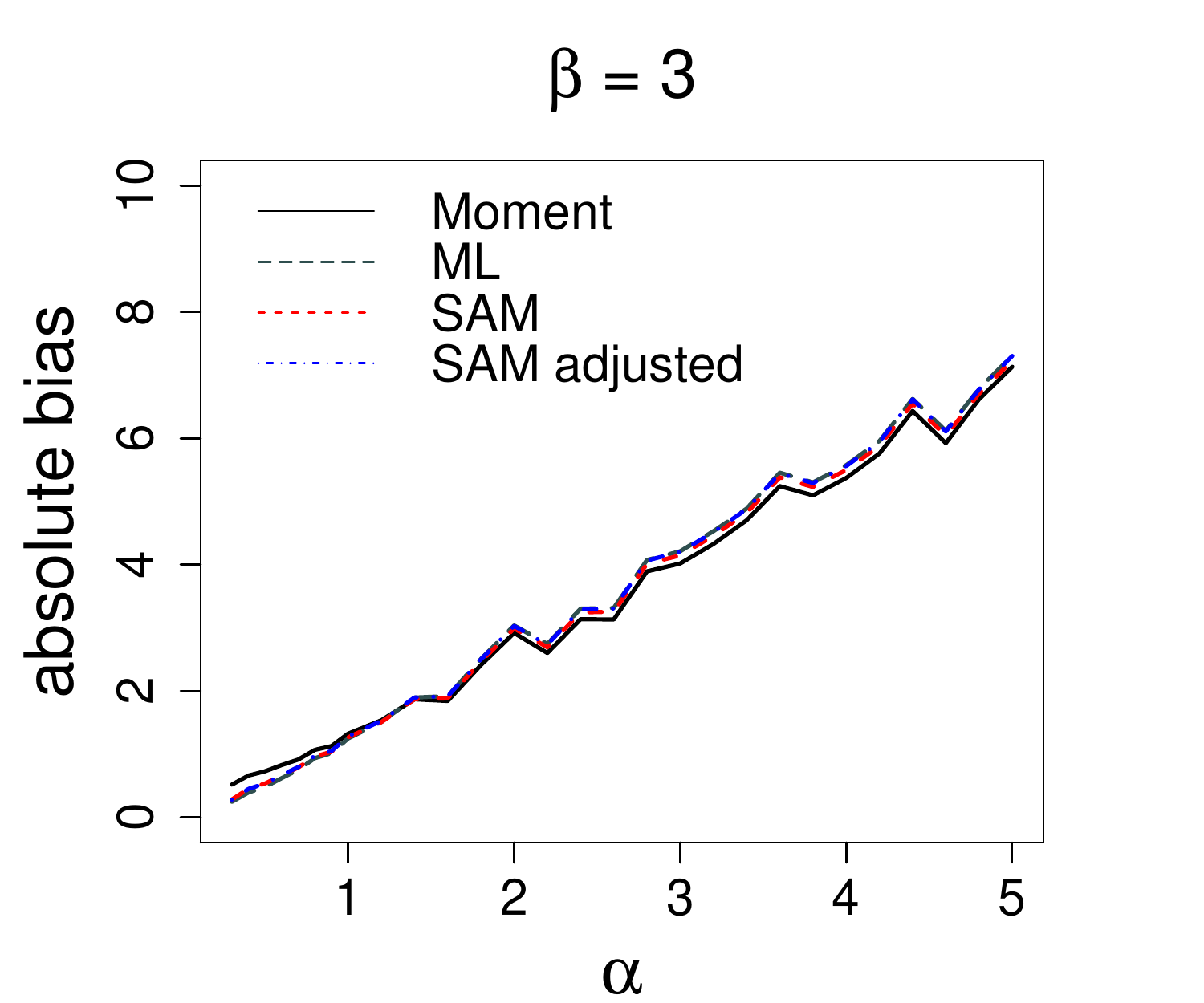}}
\subfigure{\includegraphics[width=7cm]{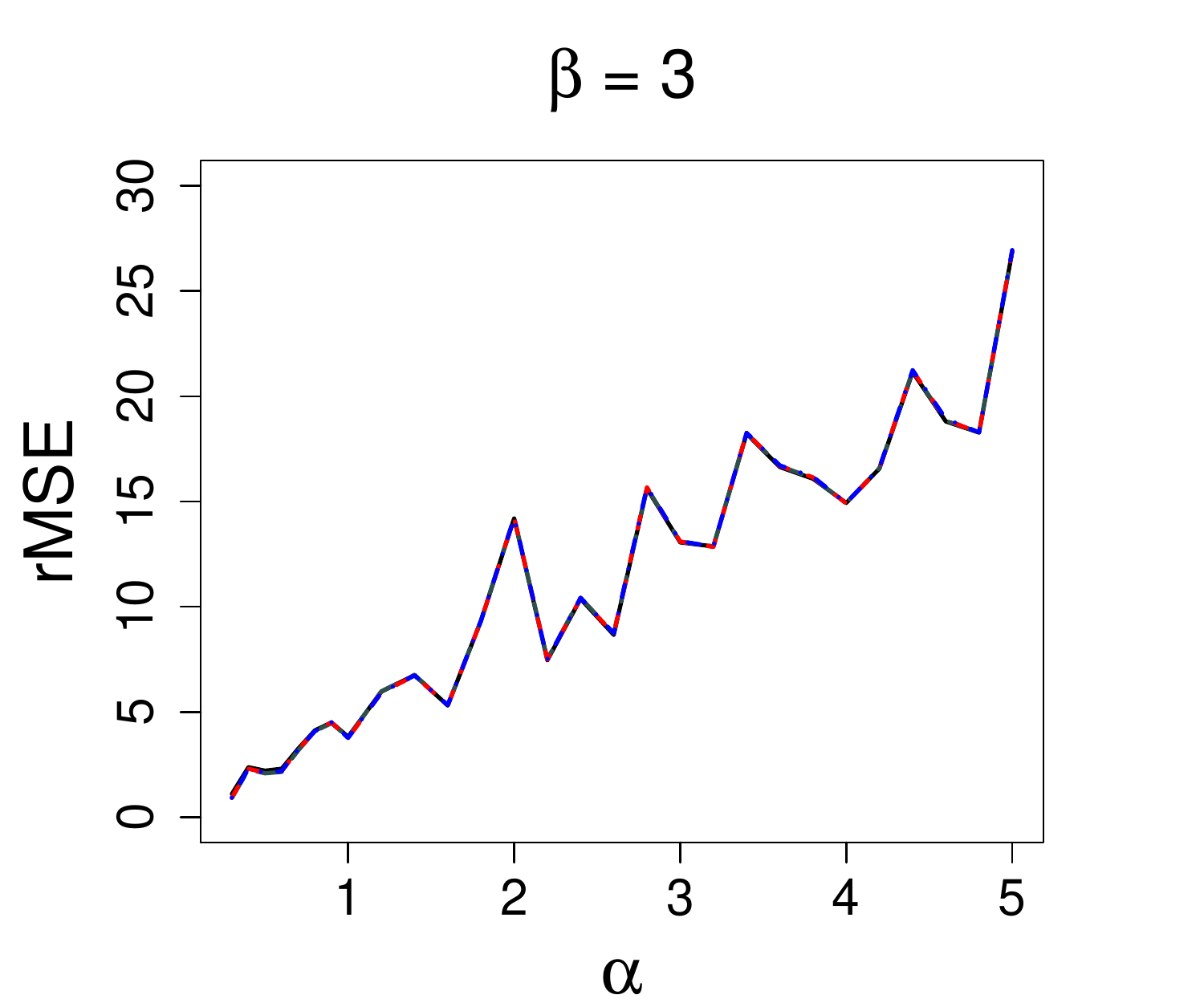}}
\caption{Absolute biases and rMSEs of the moment estimator, the ML estimator and the two proposed estimatorsunder $\beta=0.5,1,3$ and $n=5$. }
\label{fig:n5}
\end{figure}

\begin{figure}
\centering
\subfigure{\includegraphics[width=7cm]{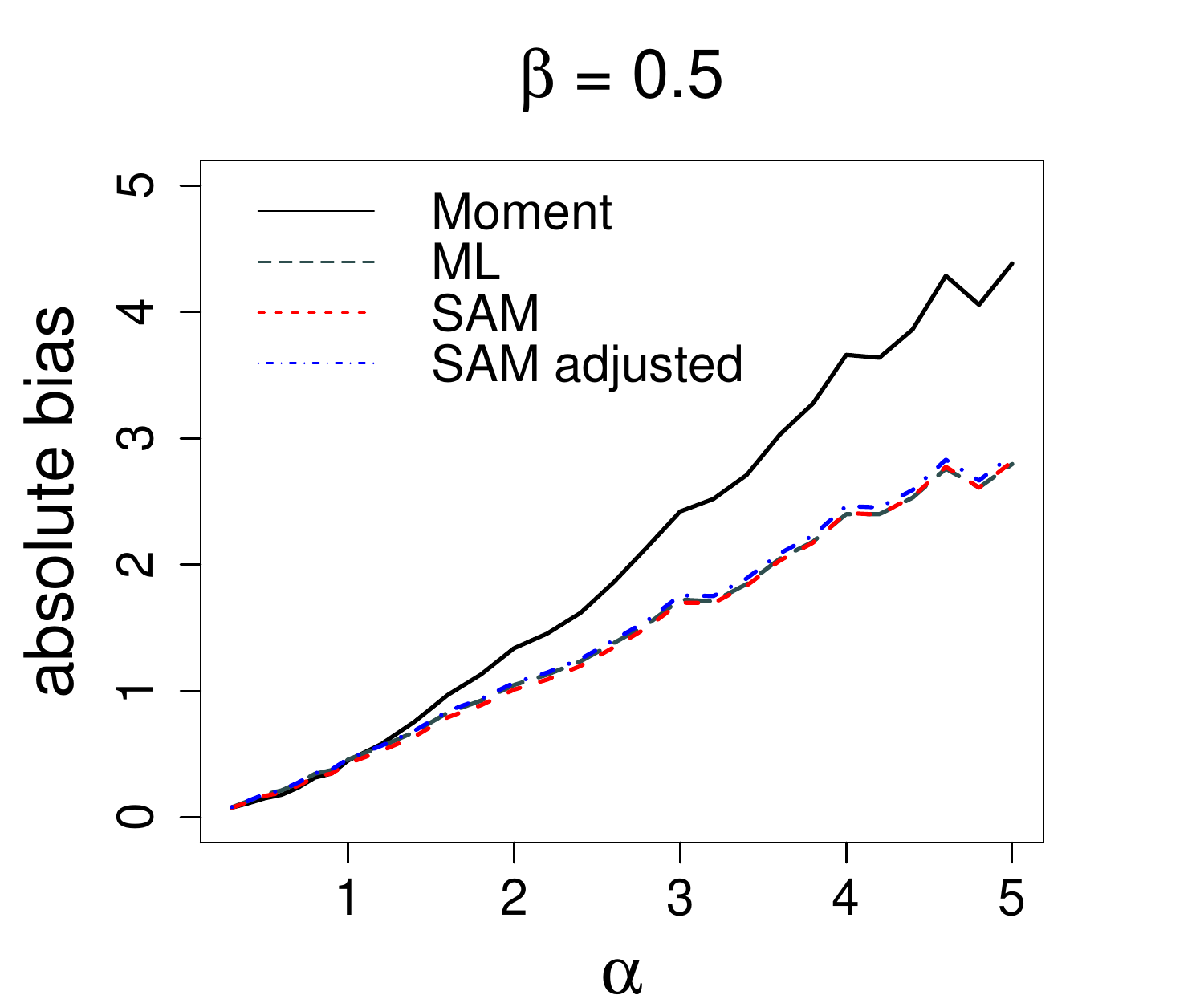}}
\subfigure{\includegraphics[width=7cm]{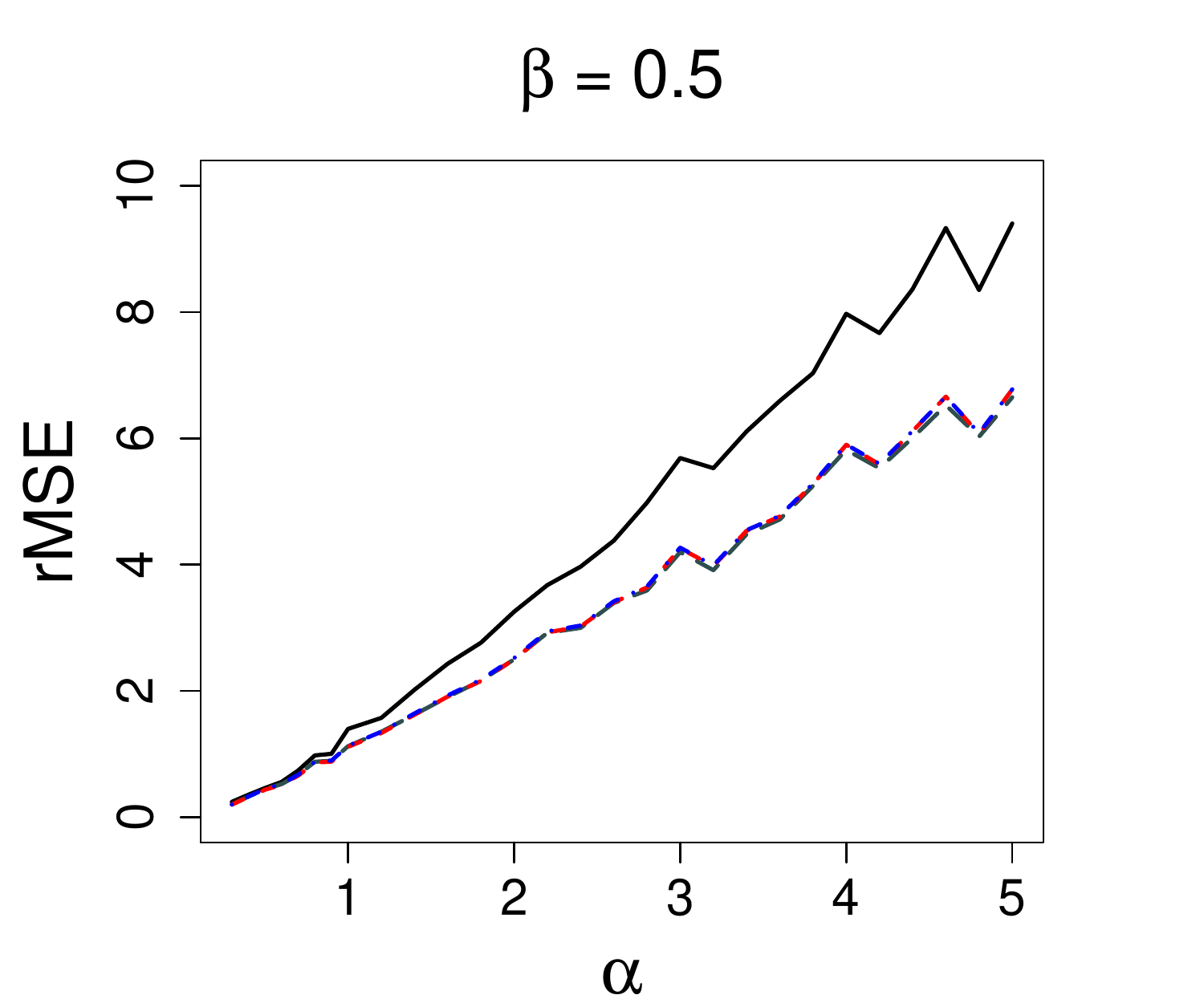}}
\hfill
\\
\subfigure{\includegraphics[width=7cm]{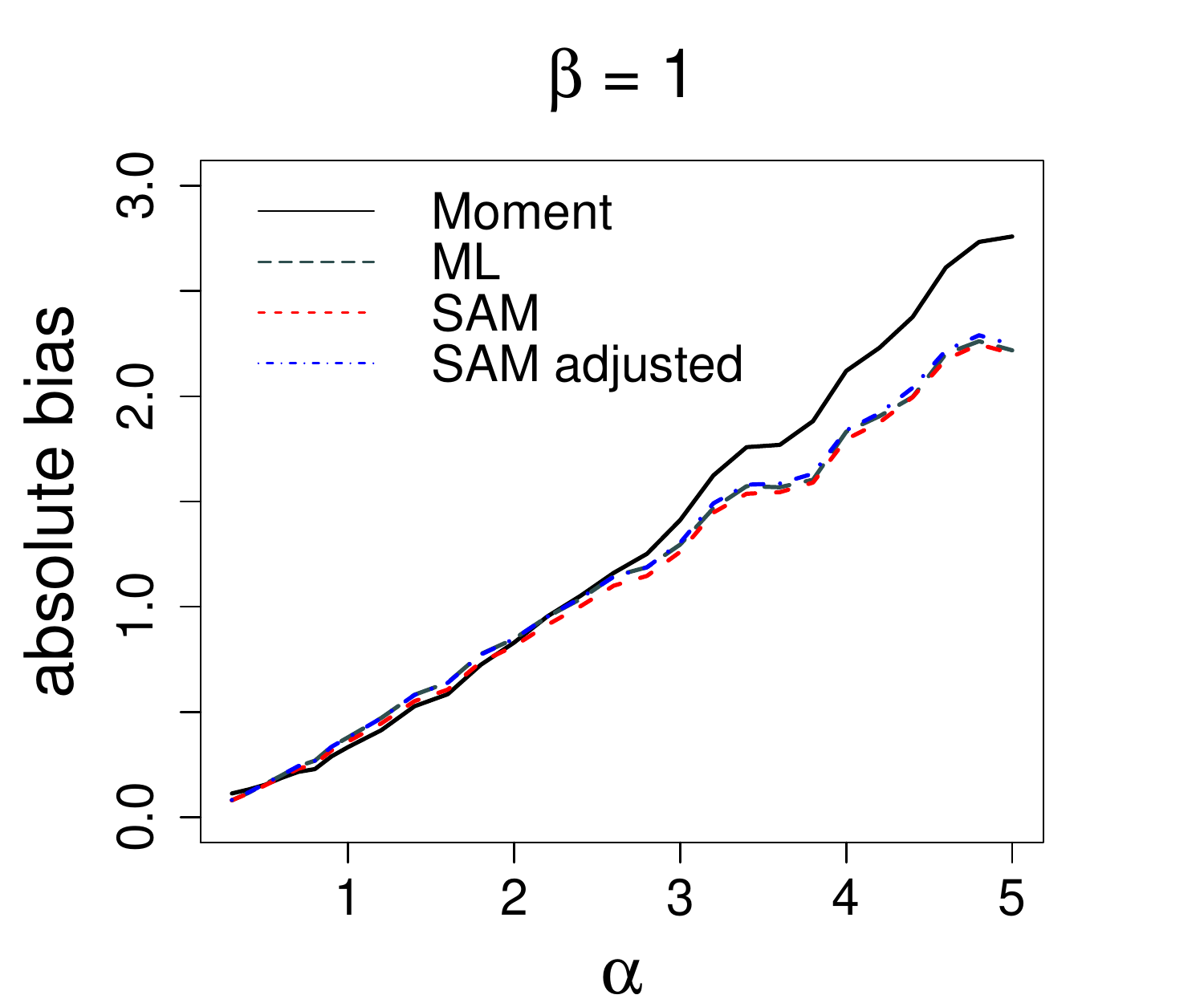}}
\subfigure{\includegraphics[width=7cm]{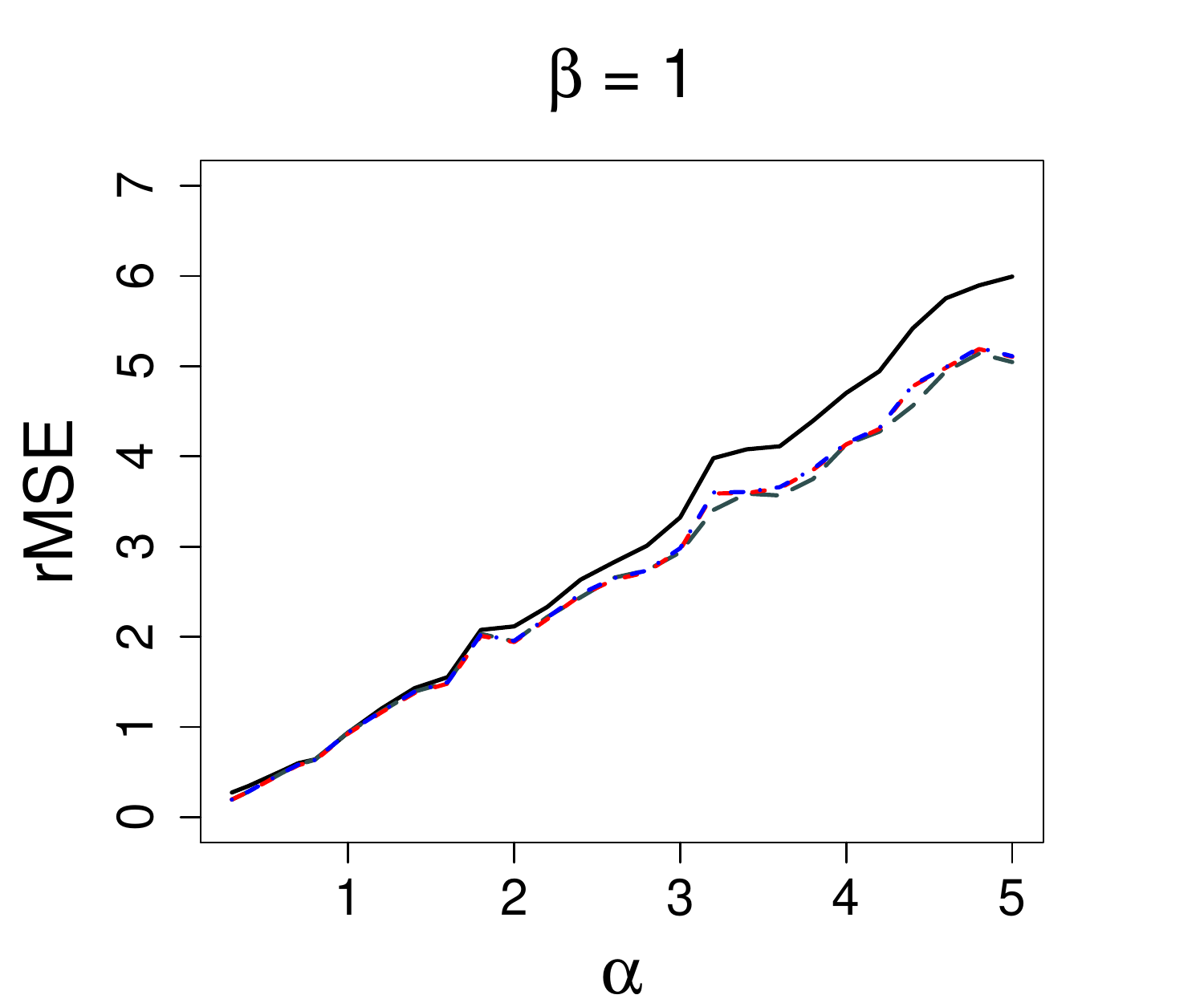}}
\hfill
\\
\subfigure{\includegraphics[width=7cm]{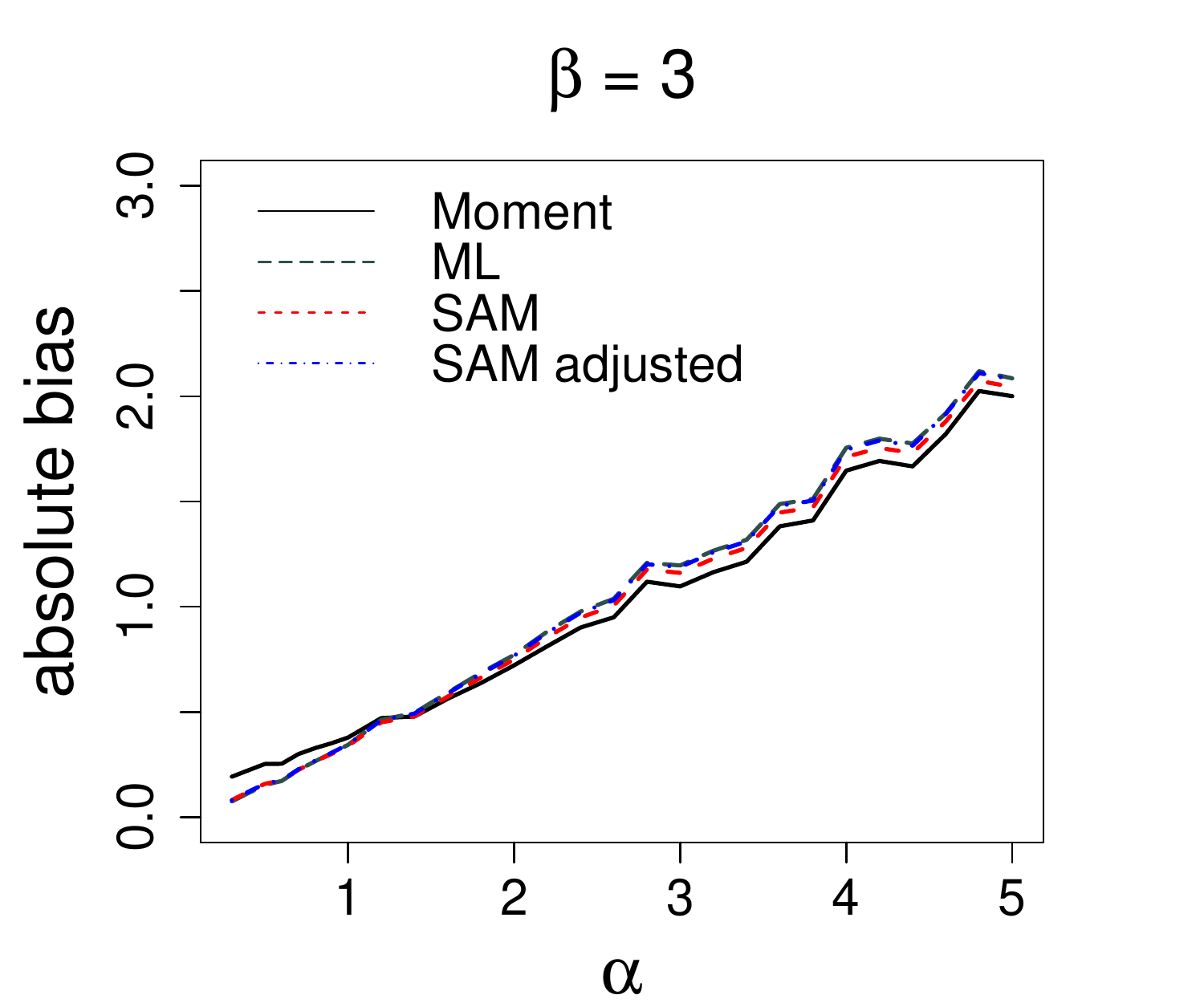}}
\subfigure{\includegraphics[width=7cm]{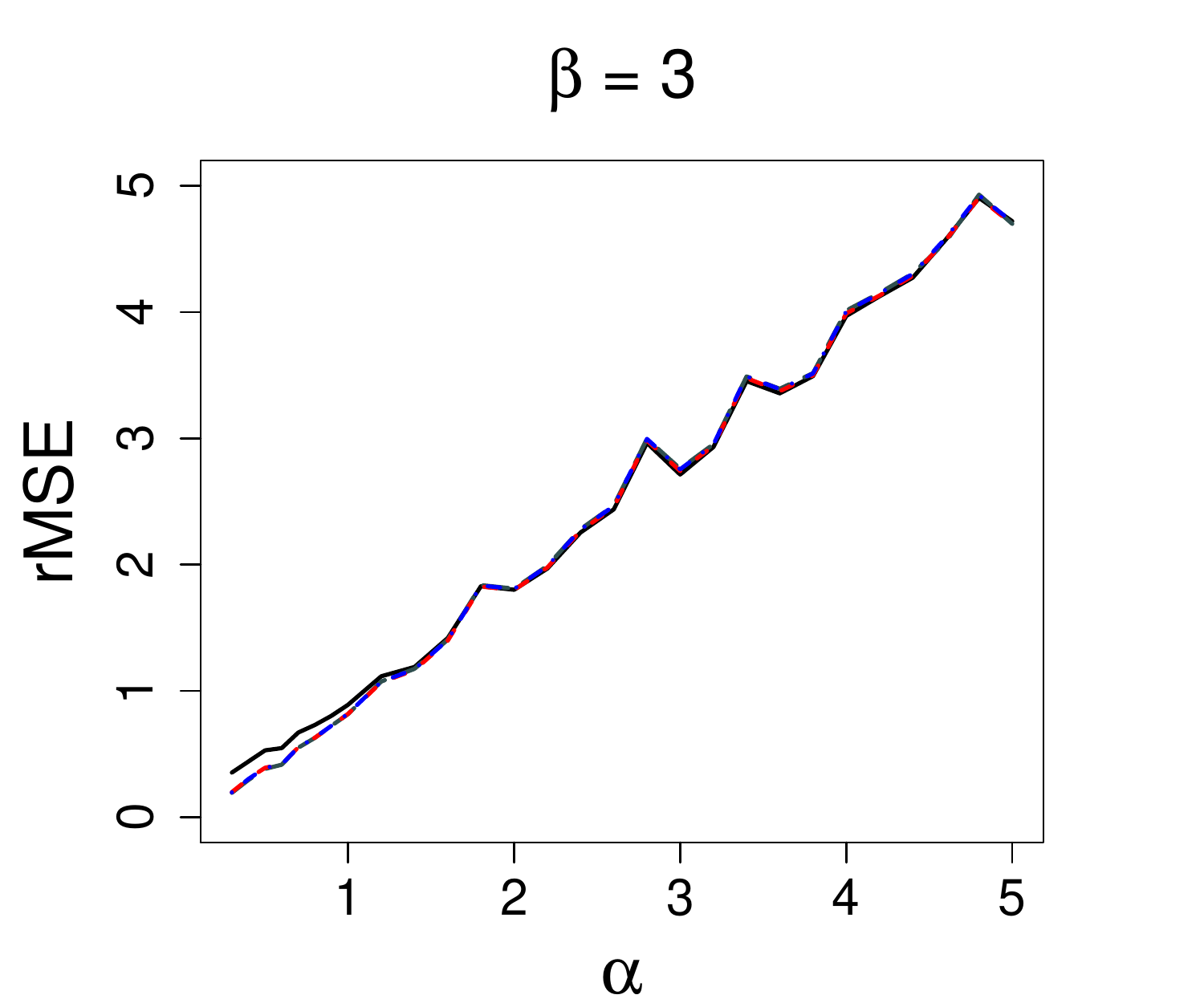}}
\caption{Absolute biases and rMSEs of the moment estimator, the ML estimator and the two proposed estimatorsunder $\beta=0.5,1,3$ and $n=10$.}
\label{fig:n10}
\end{figure}

\begin{figure}
\centering
\subfigure{\includegraphics[width=7cm]{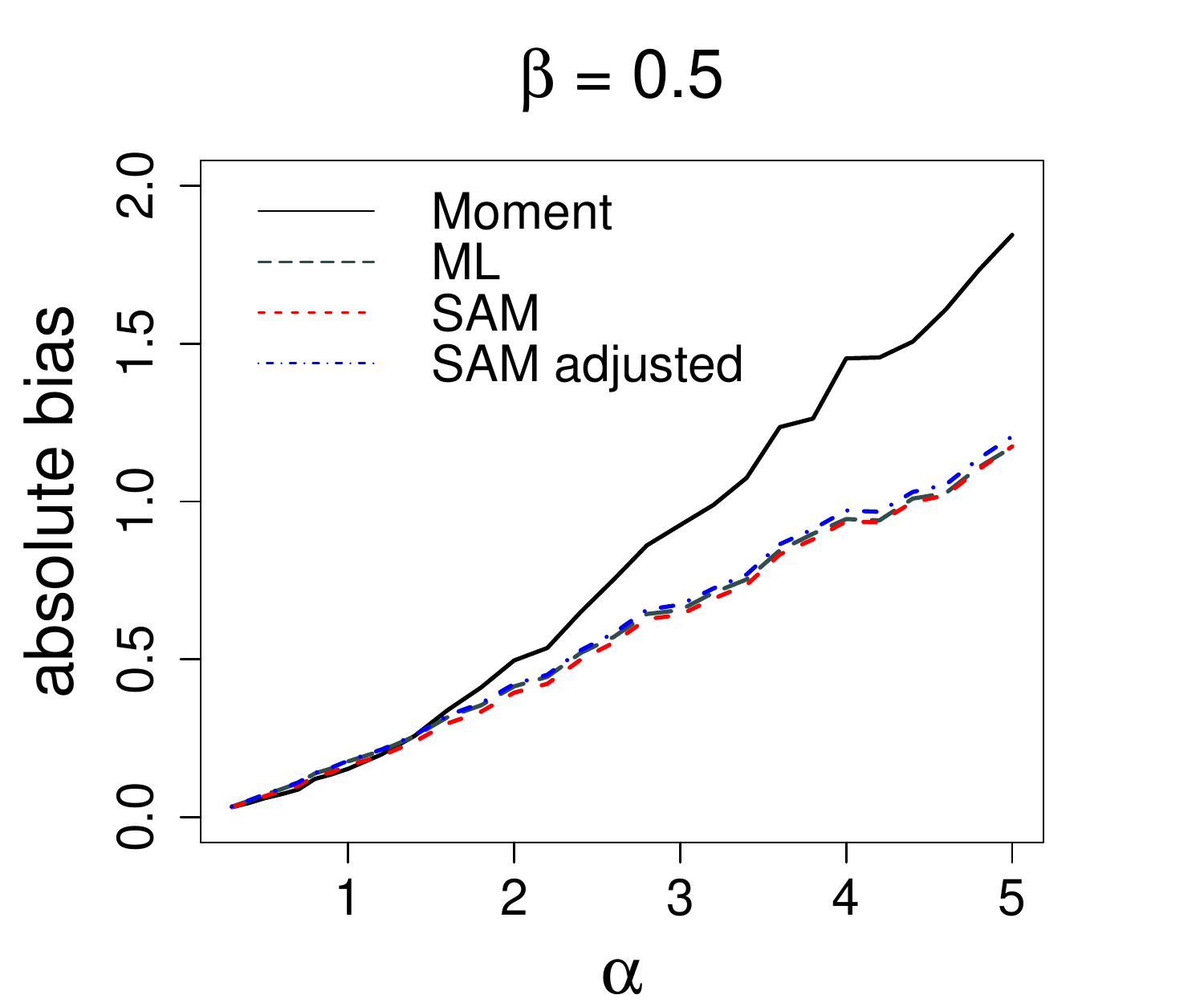}}
\subfigure{\includegraphics[width=7cm]{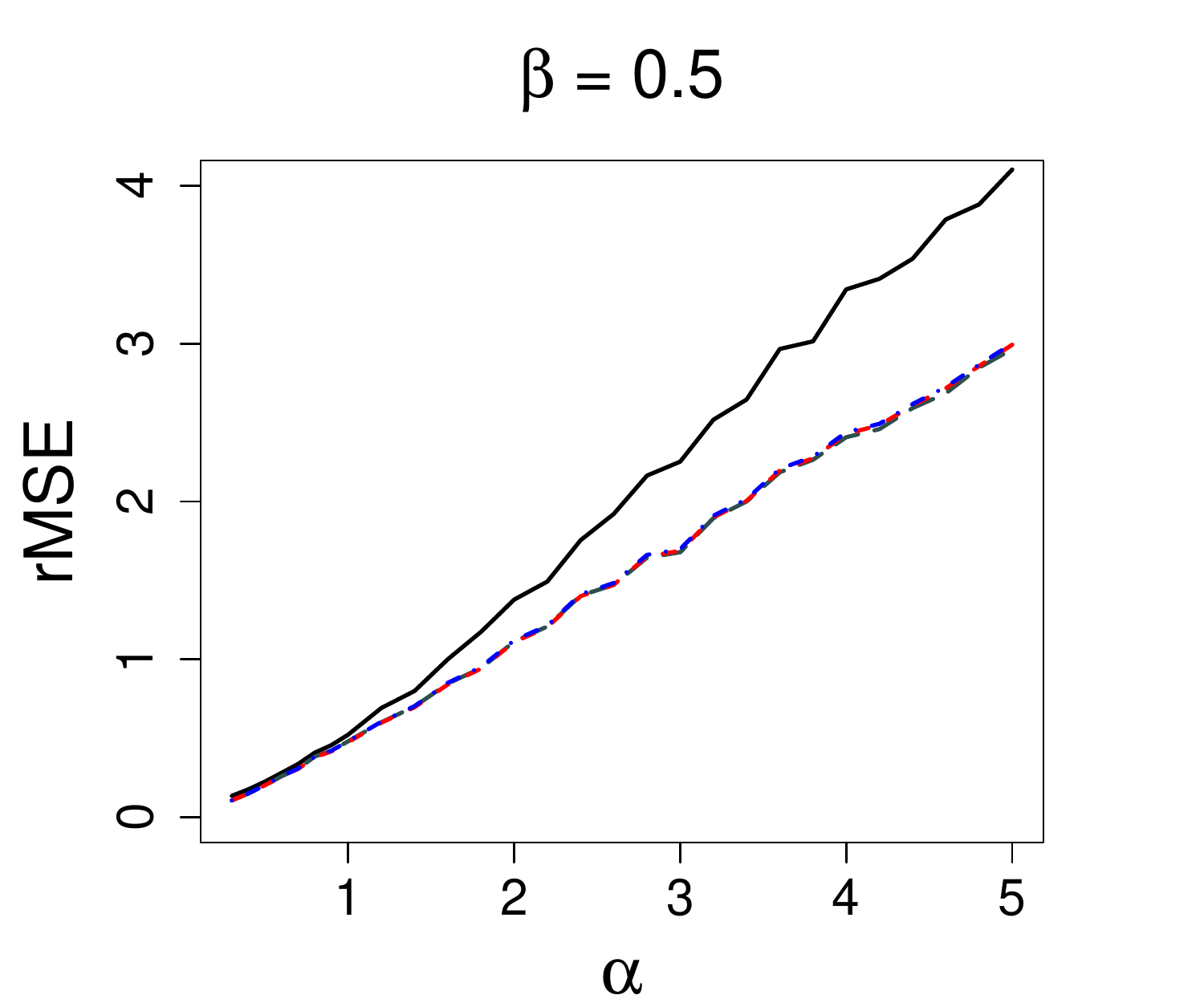}}
\hfill
\\
\subfigure{\includegraphics[width=7cm]{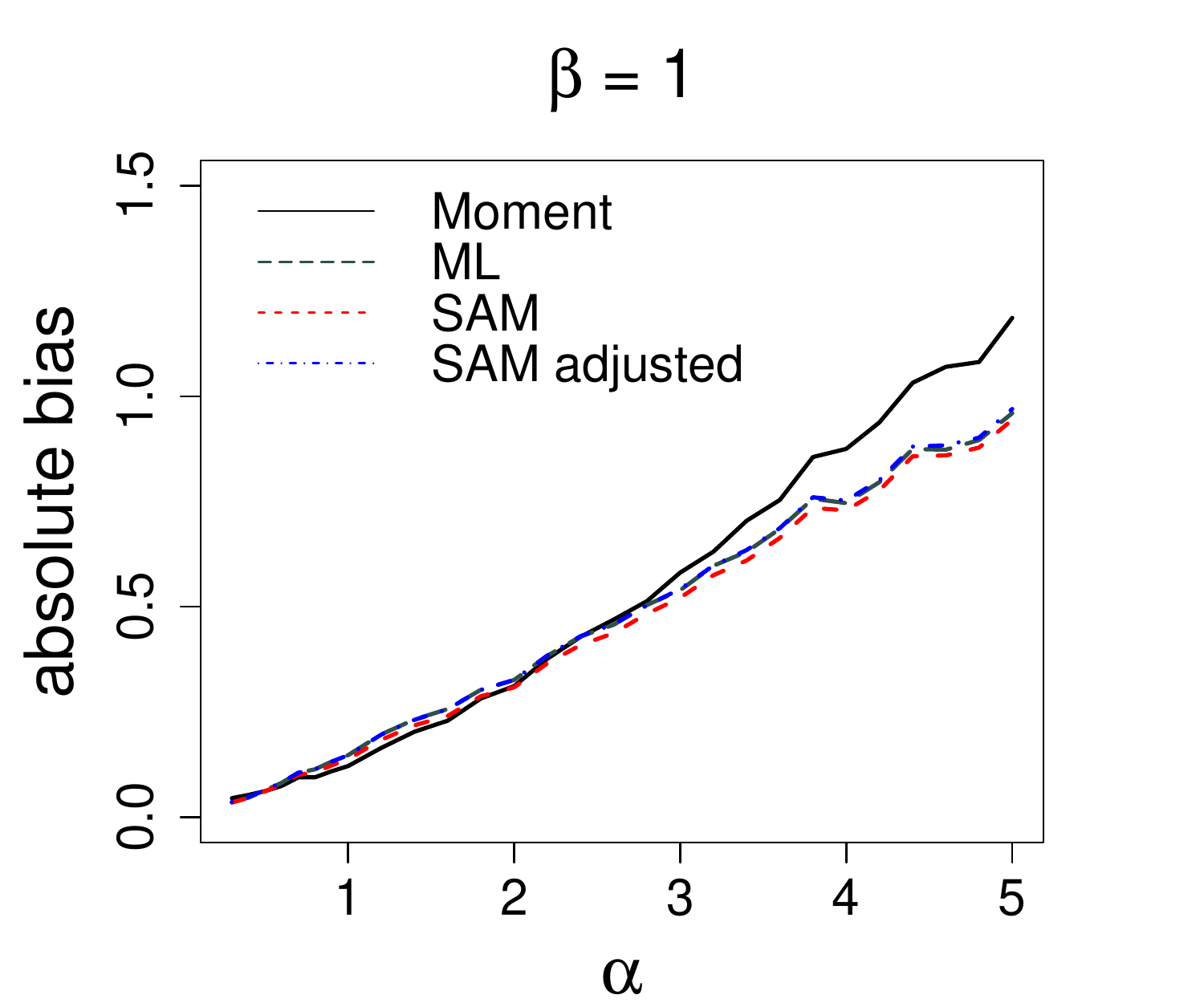}}
\subfigure{\includegraphics[width=7cm]{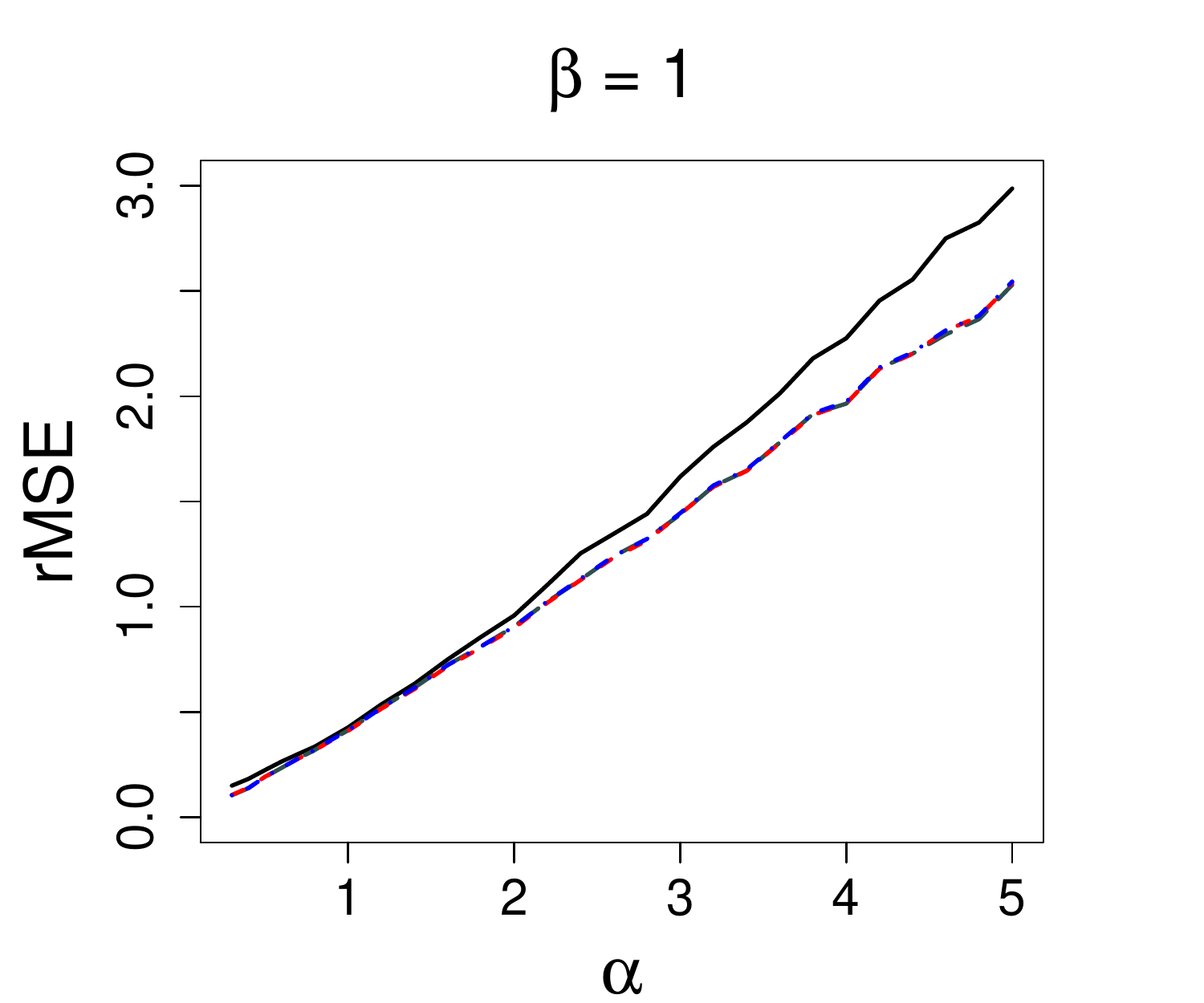}}
\hfill
\\
\subfigure{\includegraphics[width=7cm]{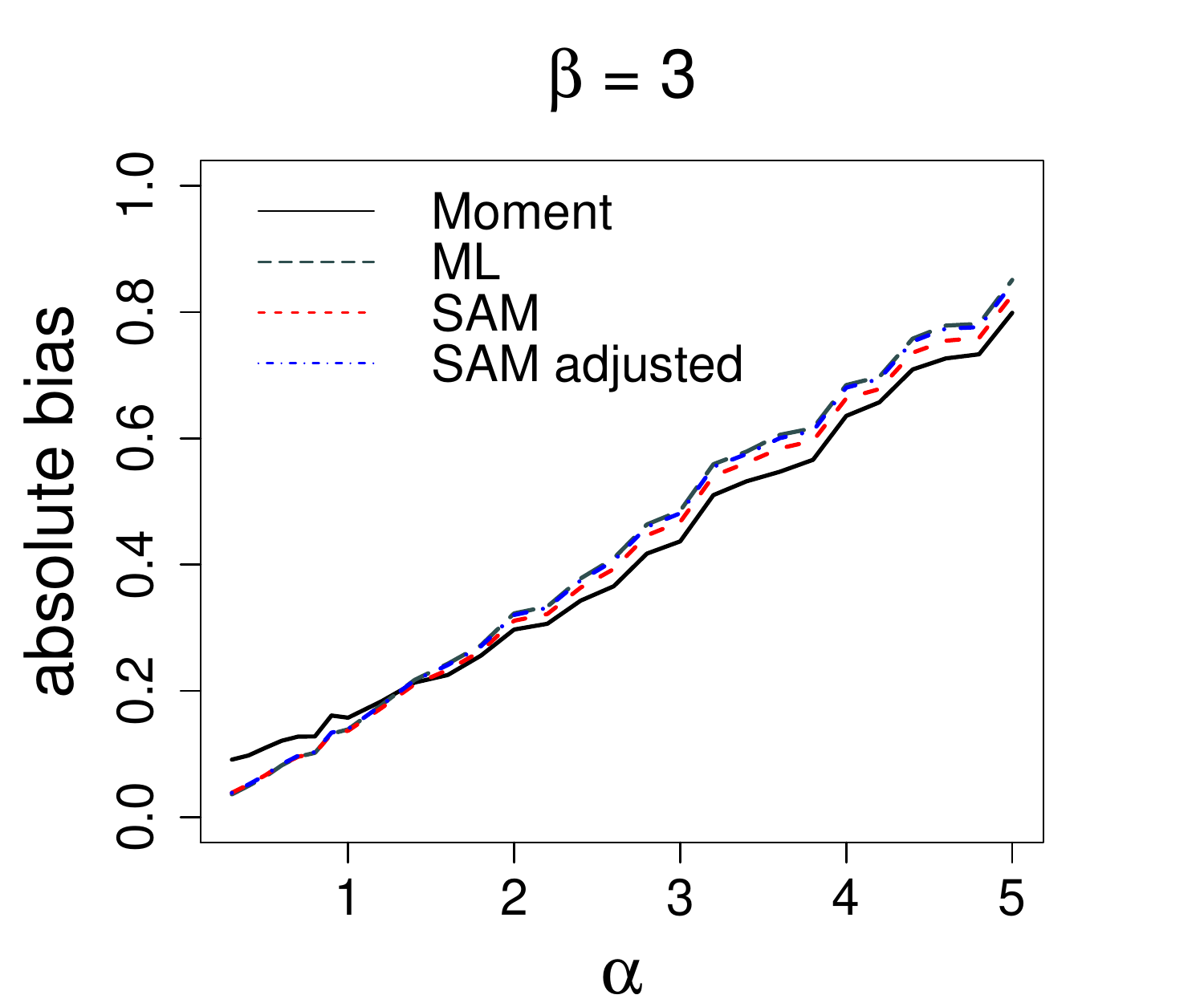}}
\subfigure{\includegraphics[width=7cm]{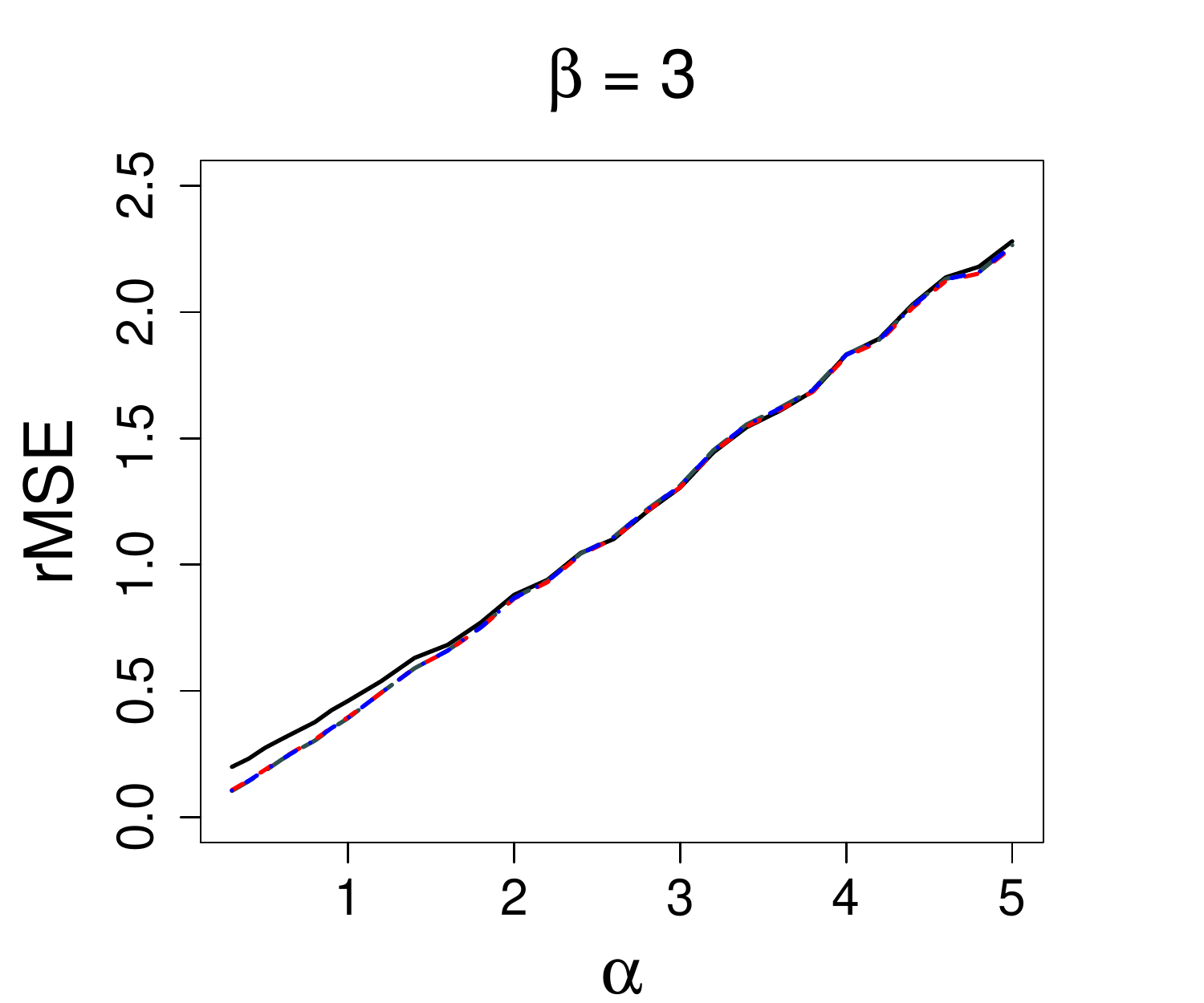}}
\caption{Absolute biases and rMSEs of the moment estimator, the ML estimator and the two proposed estimatorsunder $\beta=0.5,1,3$ and $n=20$.}
\label{fig:n20}
\end{figure}

We assess the performance of the two proposed closed-form estimators, the ML estimators and the moment estimators by simulation.  
Similar to the numerical studies in \cref{sec:idea}, we consider estimators of $\alpha$ under $\beta=0.5,1,3$. The absolute biases and the root mean square errors (rMSEs) of the four estimators are obtained based on 10000 replications, and they are shown in Figures~\ref{fig:n5}-\ref{fig:n20} under the sample sizes $n=5,10,20$, respectively. As seen, the moment estimators generally have the largest biases and MSEs among all the considered estimators. 
In addition, the two proposed estimators perform almost identically to the ML estimators in all the scenarios. It is important to mention that the quasi-Newton method is used to obtain the ML estimators of the beta distribution where the moment estimates are set as the staring values.  Throughout the simulations, it is observed that the quasi-Newton algorithm occasionally fails to converge especially in small sample sizes. For example, the ML estimators are not available in 3.7\% of the simulation cases when $n=5$ and $(\alpha,\beta)=(3,0.5)$. This is because it is very likely that $\exp(\overline{\ln X}) + \exp(\overline{\ln Y}) > 0.95$ under the setting, and the ML estimators will vary rapidly in this case \citep{cordeiro1997bias}. On the other hand, the proposed estimators have closed-forms and hence can be efficiently computed in all the cases.


\section{Conclusion}
\label{sec:conclusion}

In this paper, we have proposed two novel closed-form estimators for the beta distribution. The first version of estimators is derived by using an equation involving the sufficient statistics as the second moment equation, and they can be treated as the mixed type of moment estimators. Interestingly, they agree with the existing score-adjusted estimators for the beta distribution in the literature. On the other hand, the second new estimators are derived by solving the score equations from the generalized beta distribution, and they are shown to be a refined version of the score-adjust estimators. We compared the proposed estimators with the classical ML estimators and moment estimators theoretically and by simulation. On the one hand, compared to the moment estimators, the proposed estimators are asymptotically more efficient and have smaller MSEs in finite samples. On the other hand, the proposed estimators perform almost identically to the ML estimators in both large and finite samples. Nevertheless,  the closed-form expressions make the proposed estimators preferable in practice as the iteration algorithm in obtaining the ML estimators may either fail to converge in small samples or require nonnegligible computational time in  some applications. 
Given the superb performance of the proposed estimators, it is  possible to extend them to ease the computational difficulties and improve the estimation accuracies in estimating  the beta regression models and the beta mixture models.


\appendix
\section*{Appendix} \label{appendix}

\phantomsection
\addcontentsline{toc}{section}{Appendix}

\renewcommand{\thesubsection}{A.\arabic{subsection}}
\renewcommand{\thesubsubsection}{A.\arabic{subsection}.\arabic{subsubsection}}

\subsection{Proof of \cref{Proposition1}} \label{proof:thm1}

In this section, we sketch the proofs for the properties presented in \cref{sec:idea}. Throughout this section, let 
$\X, \LX, \LY, \XLX$ and $\XLY$ be the sample mean of $X$, $\ln X$, $\ln Y$, $X\ln X$ and $X \ln Y$, respectively. 

\subsubsection{Proof of \cref{Proposition1} (i)}  \label{proof:positivity}

First, we observe that $\X$ and $\Y$, i.e., the numerators in \eqref{eq:estalpha} and \eqref{eq:estbeta} are always well-defined and positive given any random sample from the open unit interval. Therefore, it is sufficient to show that $\overline{X \ln X}-\overline{X} \cdot \overline{\ln X}+\overline{Y \ln Y}-\overline{Y} \cdot \overline{\ln Y}\coloneqq \tgam$, i.e., the denominator in \eqref{eq:estalpha} and \eqref{eq:estbeta}, is positive. By rearranging the terms in $\tgam$, we obtain that 
\begin{align*}
\tgam&=\frac{1}{n}\left[n\XLX-n\X \cdot \LX+n\YLY-n\Y \cdot \LY\right]\\
&=\frac{1}{n}\left[\sum_{i=1}^{n}X_i\ln X_i-n\X \cdot \LX +\sum_{i=1}^{n}Y_i\ln Y_i-n\Y \cdot \LY\right]\\
&=\frac{1}{n}\left[\sum_{i=1}^{n}\left(X_i-\X\right)  \left(\ln X_i-\LX\right)+\sum_{i=1}^{n}\left(Y_i-\Y\right)\left(\ln Y_i-\LY\right)\right]\\
&=\frac{1}{n}\left[\sum_{i=1}^{n}\left(X_i-\X\right)  \left(\ln X_i-\ln \X\right)+\sum_{i=1}^{n}\left(Y_i-\Y\right)\left(\ln Y_i-\ln \Y\right)\right].
\end{align*}

Since $\ln(\cdot)$ is strictly increasing we have $(X_i - \X)(\ln{X_i} - \ln{\X}) \geq 0$ with equality if and only $X_i = \X$ for $i=1,2,\cdots,n$ and moreover $(Y_i - \Y)(\ln{Y_i} - \ln{\Y}) \geq 0$ for all $i=1,2,\cdots, n$. Therefore it follows that $\tilde{\gamma}\geq 0$ with equality if and only if $X_i=\X$ for all $i=1,2,\cdots,n$. Finally, since there exists at least one pair of distinct observations we obtain that $\tilde{\gamma} > 0$, concluding the proof.


\subsubsection{Proof of \cref{Proposition1} (ii)} \label{proof:consistency}

First of all, we define three functions 
\begin{align*}
h(a,b,c,d,e) &=d-ab-e+ac,\\
h_1(a,b,c,d,e) &= \frac{a}{h(a,b,c,d,e)},\\
h_2(a,b,c,d,e) &= \frac{1-a}{h(a,b,c,d,e)}.
\end{align*}
Then, the proposed estimators $\tilde \alpha$ and $\tilde \beta$ in \eqref{eq:estalpha} and \eqref{eq:estbeta} can be expressed as 
$$\tilde \alpha = h_1(\X, \LX, \LY, \XLX, \XLY) \quad \text{and} \quad \tilde \beta = h_2(\X, \LX, \LY, \XLX, \XLY).$$

Then, we show the expectations of $\X, \LX, \LY, \XLX$ and $\XLY$ as follows. Generally, these expectations can be obtained by using the corresponding moment generating functions, and they are given by
\begin{align*}
\E[\X] & = \frac{\alpha}{\alpha+\beta},\\
\E[\LX] & = \psi(\alpha)-\psi(\alpha+\beta),\\
\E[\LY] & = \psi(\beta)-\psi(\alpha+\beta),\\
\E[\XLX] & = \frac{\alpha}{\alpha+\beta}[\psi(\alpha+1)-\psi(\alpha+\beta+1)],\\
\E[\XLY] & = \frac{\alpha}{\alpha+\beta}[\psi(\beta)-\psi(\alpha+\beta+1)].
\end{align*}

Denote $\left(\E[\X],\E[\LX],\E[\LY], \E[\XLX],\E[\XLY]\right)^\T$ by $\bmu_1$. Based on the strong law of large numbers,  we have
$$ \left(\X, \LX, \LY, \XLX,\XLY \right)^\T \xrightarrow{a.s.} \bmu_1,$$
where $a.s.$ denotes the almost sure convergence. 
Since the denominator $h(\bmu_1)={1}/{(\alpha+\beta)}$ $>0$, the functions $h_1(\cdot)$ and $h_2(\cdot)$ are well-defined and continuous at point $\bmu_1$. An application of the continuous mapping theorem yields
$$
 \binom{\tilde \alpha}{ \tilde \beta}  \xrightarrow{a.s.} \binom{ h_1\left(\bmu_1\right)}{ h_2\left(\bmu_1\right)}.
$$
After some algebraic manipulation, we can show that $\alpha = h_1\left(\bmu_1\right)$ and $\beta = h_2\left(\bmu_1\right)$. This completes the proof of strong consistency.

\subsubsection{Proof of \cref{Proposition1} (iii)} \label{proof:normality}

Based on the central limit theorem, we have 
$$
\sqrt{n} \left[\left(\X, \LX, \LY, \XLX, \XLY\right)^\T -\bmu_1 \right] \to^{d} \mathsf{MVN}(\boldsymbol{0}, \boldsymbol \Lambda).
$$
Here, the covariance matrix $\boldsymbol  \Lambda$ is given by
\begin{equation*}
\left[\begin{array}{ccccc}
\var(X) &  & &  &  \\
\cov(X, \ln X) & \var(\ln X) &  &  & \\
\cov(X,\ln Y) & \cov(\ln X, \ln Y)& \var(\ln Y) &  & \\
\cov(X, X\ln X) & \cov(\ln X, X \ln X) & \cov(\ln Y, X \ln X)  & \var(X \ln X) &  \\
\cov(X, X \ln Y) & \cov(\ln X, X \ln Y)  & \cov(\ln Y, X \ln Y) & \cov(X \ln X, X \ln Y) & \var(X\ln Y)
 \end{array}\right]
\end{equation*}
where the upper triangular part is suppressed as it can be found easily by the symmetry.
To figure out the asymptotic covariance matrix, we need to calculate the variances and covariances involving $X$, $\ln X$, $\ln Y$, $X\ln X$, $X\ln Y$ in terms of $\alpha$ and $\beta$. To save space, the expressions are summarized in \cref{betamoments} and the detailed calculations are omitted here. 
Recall that $\tilde \alpha = h_1(\X, \LX, \LY, \XLX, \XLY)$ and $\tilde \beta = h_2(\X, \LX, \LY, \XLX, \XLY)$. The delta method therefore implies that 
\begin{equation*}
\sqrt{n} \left[\binom{\tilde \alpha}{\tilde \beta}-\binom{\alpha}{\beta}\right] \to^{d} \mathsf{MVN}(\boldsymbol{0},Q \boldsymbol \Lambda Q^\T),
\end{equation*}
where 
\begin{align*}
Q & =
\left[\begin{array}{ccccc}
\frac{\partial h_1}{\partial a} & \frac{\partial h_1}{\partial b} & \frac{\partial h_1}{\partial c} & \frac{\partial h_1}{\partial d}
& \frac{\partial h_1}{\partial e}\\
\frac{\partial h_2}{\partial a} & \frac{\partial h_2}{\partial b} & \frac{\partial h_2}{\partial c} & \frac{\partial h_2}{\partial d}
& \frac{\partial h_2}{\partial e}
 \end{array}\right] \Biggr\rvert_{(a,b,c,d,e)^\T=\bmu_1}\\
 & = \left[\begin{array}{ccccc}
\alpha(\alpha+\beta)[\psi(\alpha+1)-\psi(\beta)] & \alpha^2 & -\alpha^2 & -\alpha(\alpha+\beta) & \alpha(\alpha+\beta)\\
-\beta(\alpha+\beta)[\psi(\beta+1)-\psi(\alpha)] & \alpha\beta & -\alpha\beta & -\beta(\alpha+\beta) & \beta(\alpha+\beta)
 \end{array}\right].
\end{align*}
After tedious manipulation, it can be shown that 
\small
\begin{align*}
 &Q \boldsymbol \Lambda Q^\T =\\ 
 &\frac{1}{\alpha+\beta+1}\left[\begin{array}{cc}{\alpha^{3} \beta[\psi_1(\alpha)+ \psi_1(\beta)]+\alpha^{2}(\alpha+\beta+1)-\alpha \beta} & 
{(\beta-1) \alpha^{2}+\alpha^{2} \beta^{2}[\psi_1(\alpha)+ \psi_1(\beta)]+(\alpha-1)\beta^{2}} \\ 
{(\beta-1) \alpha^{2}+\alpha^{2} \beta^{2}[\psi_1(\alpha)+ \psi_1(\beta)]+(\alpha-1)\beta^{2}}  &
 {\alpha \beta^{3}[\psi_1(\alpha)+ \psi_1(\beta)]+\beta^{2}(\alpha+\beta+1)-\alpha \beta}\end{array}\right],
\end{align*}
which completes the proof.

\subsection{Proof of \cref{Proposition2}} \label{proof:thm2}

\subsubsection{Proof of \cref{Proposition2} (i)} \label{proof:positivity2}

We first construct two auxiliary functions on the unit interval $(0,1)$ as
\begin{align*}
g_1(x)&=\frac{x\ln x}{1-x}~~~, & g_2(x)&=\frac{x\ln x}{(1-x)\ln(1-x)}.
\end{align*}
It is easy to show that, for any $x\in(0,1)$,
\begin{align*}
g'_1(x)&=\frac{\ln x+1-x}{(1-x)^2}<0, \\
g'_2(x)&=\frac{x\ln x + (1-x+\ln x)\ln(1-x)}{(1-x)^2\ln^2(1-x)}<0,
\end{align*}
and $\lim_{x\to 0} g_1(x)=0,\lim_{x\to 1} g_1(x)=-1$.

Following the boundedness and monotonicity of $g_1$, we obtain that $-1\le M_X,M_Y\le 0$. Therefore, two numerators in \cref{eq:estalpha2,eq:estbeta2} satisfy
\begin{align*}
(1+M_X)\LY+(1+M_Y)M_X&<0,\\
(1+M_Y)\LX+(1+M_X)M_Y&<0.
\end{align*}

The rest is to show that the common denominator $M_XM_Y-\LX\cdot\LY$ is always negative which is equivalent to show that $\sum_i\ln X_i \sum_j\ln Y_j - \sum_i\frac{X_i\ln X_i}{1-X_i} \sum_j\frac{Y_j\ln Y_j}{1-Y_j} >0$. By expanding the products, we have 
\begin{align*}
 &\sum_i\ln X_i \sum_j\ln Y_j - \sum_i\frac{X_i\ln X_i}{1-X_i} \sum_j\frac{Y_j\ln Y_j}{1-Y_j}\\
=&\sum_i\sum_j\left[\ln X_i \ln (1-X_j)-\frac{X_i\ln X_i}{1-X_i}\frac{(1-X_j)\ln (1-X_j)}{X_j}\right]\\
=&\sum_i\sum_{j\ne i}\left[\frac{X_j-X_i}{(1-X_i)X_j}\ln X_i\ln (1-X_j)\right]\\
=&\sum_i\sum_{j>i}\left[\frac{X_j-X_i}{(1-X_i)X_j}\ln X_i\ln (1-X_j) + \frac{X_i-X_j}{(1-X_j)X_i}\ln X_j\ln (1-X_i) \right]\\
=&\sum_i\sum_{j>i}(X_j-X_i)\left[\frac{\ln X_i\ln (1-X_j)}{(1-X_i)X_j} - \frac{\ln X_j\ln (1-X_i)}{(1-X_j)X_i} \right]\\
=&\sum_i\sum_{j>i}(X_j-X_i)\frac{\ln(1-X_i)\ln(1-X_j)}{X_iX_j}\left[\frac{X_i\ln X_i}{(1-X_i)\ln(1-X_i)} - \frac{X_j\ln X_j}{(1-X_j)\ln(1-X_j)} \right]\\
=&\sum_i\sum_{j>i}(X_j-X_i)\frac{\ln(1-X_i)\ln(1-X_j)}{X_iX_j}\left[g_2(X_i)- g_2(X_j) \right]\ge 0,
\end{align*}
where the last inequality follows the monotonicity of $g_2$. Finally, since there exists at least one pair of distinct observations we obtain that the above inequality holds strictly, concluding the proof. 

\subsubsection{Proof of \cref{Proposition2} (ii)} \label{proof:consistency2}

First of all, we define two functions 
\begin{align*}
h_3(a,b,c,d) &= \frac{(1+c)b+(1+d)c}{cd-ab},\\
h_4(a,b,c,d) &= \frac{(1+d)a+(1+c)d}{cd-ab}.
\end{align*}
Then, the proposed estimators $\bralpha$ and $\brbeta$ in \eqref{eq:estalpha2} and \eqref{eq:estbeta2} can be expressed as 
$$\bralpha = h_3(\LX, \LY, M_X,M_Y) \quad \text{and} \quad \brbeta = h_4(\LX, \LY, M_X,M_Y).$$

We first show the expectations of $\LX, \LY, M_X$ and $M_Y$ as follows. Generally, these expectations can be obtained by using the normalization condition of beta densities, and they are given by
\begin{align*}
\E[\LX] & = \psi(\alpha)-\psi(\alpha+\beta),\\
\E[\LY] & = \psi(\beta)-\psi(\alpha+\beta),\\
\E[M_X] & = \frac{\alpha}{\beta-1}[\psi(\alpha+1)-\psi(\alpha+\beta)],\\
\E[M_Y] & = \frac{\beta}{\alpha-1}[\psi(\beta+1)-\psi(\alpha+\beta)].
\end{align*}

One shall notice that both $\E[M_X]=\alpha\psi_1(\alpha+1)$ if $\beta=1$ and $\E[M_Y]=\beta\psi_1(\beta+1)$ if $\alpha=1$ are still well-defined since they are just removable discontinuities. 

Denote $\left(\E[\LX],\E[\LY], \E[M_X] ,\E[M_Y] \right)^\T$ by $\bmu_2$. Based on the strong law of large numbers,  we have
$$ \left(\LX, \LY, M_X, M_Y \right)^\T \xrightarrow{a.s.} \bmu_2,$$
where $a.s.$ denotes the almost sure convergence. 

If the common denominator of $h_3(\bmu_2)$ and $h_4(\bmu_2)$, i.e. $\E[M_X]\E[M_Y]-\E[\LX]\E[\LY]$, is nonzero, the functions $h_3(\cdot)$ and $h_4(\cdot)$ are well-defined and continuous at point $\bmu_2$. Unlikely the proof in \cref{proof:consistency}, evaluating this denominator involves complicated digamma functions and becomes rather inconvenient. Therefore, we present a direct proof which is similar to \cref{proof:positivity2} as follows.
\begin{align*}
&\E[\LX]\E[\LY]-\E[M_X]\E[M_Y]\\
=& \int_{0}^{1}\ln x f(x)dx\int_{0}^{1}\ln (1-x) f(x)dx -  \int_{0}^{1}\frac{x\ln x}{1-x} f(x)dx\int_{0}^{1}\frac{(1-x)\ln (1-x)}{x} f(x)dx\\
=& \int_{0}^{1}\int_{0}^{1}\ln x \ln (1-y) f(x)f(y)dxdy -  \int_{0}^{1}\int_{0}^{1}\frac{x(1-y)}{(1-x)y}\ln x \ln (1-y) f(x)f(y)dxdy\\
=& \int_{0}^{1}\int_{0}^{1}\frac{y-x}{(1-x)y}\ln x \ln (1-y) f(x)f(y)dxdy \\
=& \int_{0}^{1}\int_{0}^{x}\frac{y-x}{(1-x)y}\ln x \ln (1-y) f(x)f(y)dxdy +\int_{0}^{1}\int_{x}^{1}\frac{y-x}{(1-x)y}\ln x \ln (1-y) f(x)f(y)dxdy \\
=& \int_{0}^{1}\int_{0}^{x}\frac{y-x}{(1-x)y}\ln x \ln (1-y) f(x)f(y)dxdy +\int_{0}^{1}\int_{0}^{y}\frac{y-x}{(1-x)y}\ln x \ln (1-y) f(x)f(y)dydx \\
=& \int_{0}^{1}\int_{0}^{x}\frac{y-x}{(1-x)y}\ln x \ln (1-y) f(x)f(y)dxdy +\int_{0}^{1}\int_{0}^{x}\frac{x-y}{(1-y)x}\ln y \ln (1-x) f(y)f(x)dxdy \\
=& \int_{0}^{1}\int_{0}^{x}(y-x)\frac{\ln(1-x)\ln(1-y)}{xy}[g_2(x)-g_2(y)]f(x)f(y)dxdy >0 .
\end{align*}

Finally, an application of the continuous mapping theorem yields
$$
 \binom{\tilde \alpha}{ \tilde \beta}  \xrightarrow{a.s.} \binom{ h_3\left(\bmu_2\right)}{ h_4\left(\bmu_2\right)}.
$$
After some algebraic manipulation, we can show that $\alpha = h_3\left(\bmu_2\right)$ and $\beta = h_4\left(\bmu_2\right)$. This completes the proof of strong consistency.

\subsubsection{Proof of \cref{Proposition2} (iii)} \label{proof:normality2}


The estimators $\bralpha$ and $\brbeta$ are the solutions to the following estimating equation involving $\mlx,\mly,M_X,M_Y$ as
\[
g(\alpha,\beta)=\binom{1+\alpha\ln X-(\beta-1)M_X}{1+\beta\ln Y-(\alpha-1)M_Y}=\bzero~.
\]
Following the asymptotic distribution of generalized method of moments \citep{hansen1982large}, the covariance matrix $\Sigma_2$ can be expressed as $\Sigma_2=G^{-1}\Omega (G^\T)^{-1}$ where
\[
G=\E\begin{bmatrix} \ln X & -M_X \\ -M_{Y} &  \ln Y\end{bmatrix},
\]
and
{\footnotesize
\[
\Omega=\E\begin{bmatrix} (1+\alpha\ln X-(\beta-1)M_X)^2 & [1+\alpha\ln X-(\beta-1)M_X][1+\beta \ln Y-(\alpha-1)M_{Y}] \\ [1+\alpha\ln X-(\beta-1)M_X][1+\beta\ln Y-(\alpha-1)M_{Y}] & (1+\beta \ln Y-(\alpha-1)M_{Y})^2 \end{bmatrix}.
\]
}

One can easily show that 
\[
G=
\begin{bmatrix} \psi(\alpha)-\psi(\alpha+\beta) & -\frac{\alpha}{\beta-1}[\psi(\alpha+1)-\psi(\alpha+\beta)] \\ 
-\frac{\beta}{\alpha-1}[\psi(\beta+1)-\psi(\alpha+\beta)] & \psi(\beta)-\psi(\alpha+\beta)  \end{bmatrix} = 
\begin{bmatrix} \kappa_{\alpha\beta} & -\tau_{\alpha\beta} \\ -\tau_{\beta\alpha}  & \kappa_{\beta\alpha}  \end{bmatrix}.
\]
For the matrix $\Omega$, two terms in the first row are given by 
{\footnotesize
\begin{align*}
&\E[(1+\alpha\ln X-(\beta-1)M_X)^2]=\var(1+\alpha\ln X-(\beta-1)M_X)\\
=&~\alpha^2\var(\ln X)+(\beta-1)^2\var(M_X)-2\alpha(\beta-1)\cov(\ln X,M_X)\\
=&~1+\frac{\alpha}{\beta-2}[(\alpha+\beta-1)(\psi_1(\alpha+1)-\psi_1(\alpha+\beta))+2(\psi(\alpha+1)-\psi(\alpha+\beta))+(\alpha+\beta-1)(\psi(\alpha+1)-\psi(\alpha+\beta))^2]\\
=&~\frac{\beta}{\beta-2}+\frac{\alpha(\alpha+\beta-1)}{\beta-2}[\kappa_1(\alpha,\beta)+\kappa_{\alpha\beta}^2]+\frac{2(2\alpha+\beta-1)}{\beta-2}\kappa_{\alpha\beta}=\omega_{\alpha\beta},\\
&\E\left[(1+\alpha\ln X-(\beta-1)M_X)(1+\beta \ln Y-(\alpha-1)M_{Y})\right]=\cov(\alpha\ln X-(\beta-1)M_X,\beta \ln Y-(\alpha-1)M_{Y})\\
=&~\alpha\beta\cov(\ln X, \ln Y)-\alpha(\alpha-1)\cov(\ln X,M_{Y})-\beta(\beta-1)\cov( \ln Y,M_X)+(\alpha-1)(\beta-1)\cov(M_X,M_{Y})\\
=&~\frac{\alpha\beta-1}{(\alpha-1)(\beta-1)}+\frac{\alpha}{\beta-1}[\psi(\alpha)-\psi(\alpha+\beta)]+\frac{\beta}{\alpha-1}[\psi(\beta)-\psi(\alpha+\beta)]\\
&+(\alpha+\beta+1)[\psi_1(\alpha+\beta)-[\psi(\alpha)-\psi(\alpha+\beta)][\psi(\beta)-\psi(\alpha+\beta)]]\\
=&~1+\frac{\alpha}{\beta-1}\kappa(\alpha+1,\beta-1)+\frac{\beta}{\alpha-1}\kappa(\beta+1,\alpha-1)+(\alpha+\beta+1)[\psi_1(\alpha+\beta)-\kappa_{\alpha\beta}\kappa_{\beta\alpha}]=\rho.
\end{align*}
}
By the symmetry of beta distribution, 
\[
\Omega= \begin{bmatrix}\omega_{\alpha\beta} & \rho \\  
\rho &  \omega_{\beta\alpha} \end{bmatrix}.
\]
Substituting $G$ and $\Omega$ into $\Sigma_2=G^{-1}\Omega (G^\T)^{-1}$ and rearranging the terms in the final matrix product yields the form of $\Sigma_2$ in \cref{Proposition2}.

\subsection{Moments of beta distribution} \label{betamoments}

Here we show the formulas that are used to obtain the asymptotic covariance matrix of the proposed estimators. In general, the following variances and covariances are derived using their definitions involving the expectations, and the expectations are obtained by using the normalization condition of beta densities,.
\begin{align*}
\var(X) & =  \frac{\alpha \beta}{(\alpha+\beta)^2(\alpha+\beta+1)},\\
\cov(X,\ln X) & = \frac{\beta}{(\alpha+\beta)^2},\\
 \cov(X, \ln Y) & = -\frac{\alpha}{(\alpha+\beta)^2}, \\
\cov(X, X\ln X) &  = \frac{\alpha(\alpha+1)}{(\alpha+\beta)(\alpha+\beta+1)}[\psi(\alpha+2)-\psi(\alpha+\beta+2)]\\
& \quad -\frac{\alpha^2}{(\alpha+\beta)^2}[\psi(\alpha+1)-\psi(\alpha+\beta+1)],\\
\cov(X, X\ln Y) & =  \frac{\alpha(\alpha+1)}{(\alpha+\beta)(\alpha+\beta+1)}[\psi(\beta)-\psi(\alpha+\beta+2)]-\frac{\alpha^2}{(\alpha+\beta)^2}[\psi(\beta)-\psi(\alpha+\beta+1)],\\
\var(\ln X) & =  \psi_1(\alpha)-\psi_1(\alpha+\beta), \\
 \var(\ln Y) & = \psi_1(\beta)-\psi_1(\alpha+\beta), \\
\cov(\ln X, \ln Y) & = -\psi_1(\alpha+\beta),\\
\cov(\ln X, X\ln X) &=  \frac{\beta}{(\alpha+\beta)^2}[\psi(\alpha+1)-\psi(\alpha+\beta+1)]+\frac{\alpha}{\alpha+\beta}[\psi_1(\alpha+1)-\psi_1(\alpha+\beta+1)],\\
\cov(\ln X, X\ln Y) &=  \frac{\beta}{\alpha+\beta}\psi_1(\alpha+\beta+1)-\frac{\beta}{(\alpha+\beta)^2}[\psi(\alpha+\beta+1)-\psi(\beta+1)]-\psi_1(\alpha+\beta), \\
\cov(\ln Y, X\ln X) &=  -\frac{\alpha}{(\alpha+\beta)^2}[\psi(\alpha+1)-\psi(\alpha+\beta+1)]-\frac{\alpha}{\alpha+\beta}\psi_1(\alpha+\beta+1), \\
\cov(\ln Y, X\ln Y) &=  \psi_1(\beta)-\psi_1(\alpha+\beta)-\frac{\beta}{\alpha+\beta}[\psi_1(\beta+1)-\psi_1(\alpha+\beta+1)]\\
& \quad -\frac{\alpha}{(\alpha+\beta)^2}[\psi(\beta+1)-\psi(\alpha+\beta+1)], \\
\var(X \ln X) &= \frac{\alpha(\alpha+1)}{(\alpha+\beta+1)(\alpha+\beta)}\left\{ [\psi(\alpha+2)-\psi(\alpha+\beta+2)]^2+\psi_1(\alpha+2)-\psi_1(\alpha+\beta+2) \right\}\\
& \quad -\frac{\alpha^2}{(\alpha+\beta)^2}[\psi(\alpha+1)-\psi(\alpha+\beta+1)]^2,\\
\var(X \ln Y) &= \frac{\alpha(\alpha+1)}{(\alpha+\beta+1)(\alpha+\beta)}\left\{ [\psi(\beta)-\psi(\alpha+\beta+2)]^2+\psi_1(\beta)-\psi_1(\alpha+\beta+2) \right\}\\
& \quad -\frac{\alpha^2}{(\alpha+\beta)^2}[\psi(\beta)-\psi(\alpha+\beta+1)]^2,\\
\cov(X\ln X, X\ln Y) &= -\frac{\alpha^2}{(\alpha+\beta)^2} [\psi(\alpha+1)-\psi(\alpha+\beta+1)][\psi(\beta)-\psi(\alpha+\beta+1)]\\
& \quad +\frac{\alpha(\alpha+1)}{(\alpha+\beta)(\alpha+\beta+1)}\left\{  [\psi(\alpha+2)-\psi(\alpha+\beta+2)][\psi(\beta)-\psi(\alpha+\beta+2)] \right.\\
& \qquad \qquad \qquad \qquad \qquad \quad  \left. -\psi_1(\alpha+\beta+2) \right\},\\
\var\left(\frac{X\ln X}{1-X}\right) &= \frac{\alpha(\alpha+1)}{(\beta-1)(\beta-2)}\left\{\psi_1(\alpha+2)-\psi_1(\alpha+\beta)+[\psi(\alpha+2)-\psi(\alpha+\beta)]^2\right\}\\
& \quad -\frac{\alpha^2}{(\beta-1)^2}[\psi(\alpha+1)-\psi(\alpha+\beta)]^2,\\
\var\left(\frac{Y\ln Y}{1-Y}\right) &= \frac{\beta(\beta+1)}{(\alpha-1)(\alpha-2)}\left\{\psi_1(\beta+2)-\psi_1(\alpha+\beta)+[\psi(\beta+2)-\psi(\alpha+\beta)]^2\right\}\\
& \quad -\frac{\beta^2}{(\alpha-1)^2}[\psi(\beta+1)-\psi(\alpha+\beta)]^2,\\
\cov\left(\frac{X\ln X}{1-X},\frac{Y\ln Y}{1-Y}\right) &= -\psi_1(\alpha+\beta)+[\psi(\alpha)-\psi(\alpha+\beta)][\psi(\beta)-\psi(\alpha+\beta)]\\
& \quad -\frac{\alpha\beta}{(\alpha-1)(\beta-1)}[\psi(\alpha+1)-\psi(\alpha+\beta)][\psi(\beta+1)-\psi(\alpha+\beta)] ,\\
\cov\left(\ln X,\frac{X\ln X}{1-X}\right) &= \frac{\alpha}{\beta-1}\left[\psi_1(\alpha+1)-\psi_1(\alpha+\beta)\right]+\frac{1}{\beta-1}\left[\psi(\alpha+1)-\psi(\alpha+\beta)\right],\\
\cov\left(\ln X,\frac{Y\ln Y}{1-Y}\right) &= -\frac{\beta}{\alpha-1}\psi_1(\alpha+\beta)-\frac{\beta}{(\alpha-1)^2}[\psi(\beta+1)-\psi(\alpha+\beta)],\\
\cov\left(\ln Y, \frac{X\ln X}{1-X}\right) &=-\frac{\alpha}{\beta-1}\psi_1(\alpha+\beta)-\frac{\alpha}{(\beta-1)^2}[\psi(\alpha+1)-\psi(\alpha+\beta)],\\
\cov\left(\ln Y,\frac{Y\ln Y}{1-Y}\right) &= \frac{\beta}{\alpha-1}\left[\psi_1(\beta+1)-\psi_1(\alpha+\beta)\right]+\frac{1}{\alpha-1}\left[\psi(\beta+1)-\psi(\alpha+\beta)\right].
\end{align*}

\bibpunct{(}{)}{;}{a}{}{,}
\bibliographystyle{asa}
\bibliography{ref}

\end{document}